\documentclass[a4paper,11pt,twoside]{amsart}
\usepackage{amsmath}
\usepackage{amsthm}
\usepackage{amssymb}
\usepackage[all]{xy}

\newtheorem{thm}{Theorem}[section]
\newtheorem{prop}[thm]{Proposition}
\newtheorem{lem}[thm]{Lemma}
\newtheorem{cor}[thm]{Corollary}

\numberwithin{equation}{section}

\theoremstyle{definition}
\newtheorem{definition}[thm]{Definition}
\newtheorem{remark}[thm]{Remark}

\renewcommand{\to}{\xymatrix@1@=15pt{\ar[r]&}}
\renewcommand{\rightarrow}{\xymatrix@1@=15pt{\ar[r]&}}
\renewcommand{\mapsto}{\xymatrix@1@=15pt{\ar@{|->}[r]&}}
\renewcommand{\twoheadrightarrow}{\xymatrix@1@=15pt{\ar@{->>}[r]&}}
\renewcommand{\hookrightarrow}{\xymatrix@1@=15pt{\ar@{^(->}[r]&}}
\newcommand{\congpf}{\xymatrix@1@=15pt{\ar[r]^-\sim&}}
\renewcommand{\cong}{\simeq}

\newcommand{\qqed}{\hspace*{\fill}$\Box$}

\newcommand{\Db}{{\rm D}^{\rm b}}

\newcommand{\Aut}{{\rm Aut}}

\newcommand{\Br}{{\rm Br}}
\newcommand{\NS}{{\rm NS}}
\newcommand{\Pic}{{\rm Pic}}

\newcommand{\rk}{{\rm rk}}
\newcommand{\coh}{{\rm\bf{Coh}}}

\newcommand{\End}{{\rm End}}
\newcommand{\Hom}{{\rm Hom}}

\newcommand{\Stab}{{\rm Stab}}
\newcommand{\iso}{\cong}

\newcommand{\dual}{^{\vee}}

\newcommand{\mor}[1][]{\xrightarrow{#1}}

\newcommand{\cat}[1]{\begin{bf}#1\end{bf}}
\newcommand{\Ext}{{\rm Ext}}

\newcommand{\Mod}[1]{{\ko_{#1}\text{-}\cat{Mod}}}

\newcommand{\cal}{\mathcal}
\newcommand{\ka}{{\cal A}}

\newcommand{\ke}{{\cal E}}
\newcommand{\kf}{{\cal F}}

\newcommand{\kh}{{\cal H}}

\newcommand{\kn}{{\cal N}}

\newcommand{\ko}{{\cal O}}
\newcommand{\kp}{{\cal P}}

\newcommand{\ks}{{\cal S}}
\newcommand{\kt}{{\cal T}}

\newcommand{\ZZ}{\mathbb{Z}}
\newcommand{\QQ}{\mathbb{Q}}
\newcommand{\RR}{\mathbb{R}}
\newcommand{\CC}{\mathbb{C}}

\newcommand{\HH}{\mathbb{H}}
\newcommand{\PP}{\mathbb{P}}

\def\ss{\sigma}

\begin{document}
\title{Stability conditions for generic K3 categories}

\author[D.\ Huybrechts, E.\ Macr\`i, and P.\ Stellari]{Daniel Huybrechts, Emanuele Macr\`i, and Paolo Stellari}

%\author{D.\ Huybrechts, E.\ Macr\`i, and P.\ Stellari}

\address{D.H.: Mathematisches Institut,
Universit{\"a}t Bonn, Beringstr.\ 1, 53115 Bonn, Germany}
\email{huybrech@math.uni-bonn.de}

\address{E.M.: Max Planck Institute for Mathematics, Vivatsgasse 7, 53111 Bonn, Germany}\email{macri@mpim-bonn.mpg.de}

\address{P.S.: Dipartimento di Matematica ``F. Enriques'',
Universit{\`a} degli Studi di Milano, Via Cesare Saldini 50, 20133
Milano, Italy} \email{stellari@mat.unimi.it}



\begin{abstract} A K3 category is by definition a Calabi--Yau category of dimension two.
Geo\-metrical\-ly K3 categories occur as bounded derived
categories of (twisted) coherent sheaves on K3 or abelian
surfaces. A K3 category is generic if there are no spherical
objects (or just one up to shift). We study stability conditions
on K3 categories as introduced by Bridgeland and prove his
conjecture about the topology of the stability manifold and the
autoequivalences group for generic twisted projective K3, abelian
surfaces, and K3 surfaces with trivial Picard group.\end{abstract}

 \maketitle

\section{Introduction}

A K3 category is by definition a Calabi--Yau category $\cat{T}$ of
dimension two, i.e.\ a $k$-linear triangulated category with
functorial isomorphisms $\Hom(E,F)\cong\Hom(F,E[2])^*$ for all
objects $E,F\in\cat{T}$.

Examples of K3 categories are provided by the bounded derived
category $\Db(X)$ of the abelian category $\coh(X)$ of all
coherent sheaves on a K3 or abelian surface $X$. Twisted
analogues  can be obtained by considering the bounded derived
category $\Db(X,\alpha)$ of the abelian category $\coh(X,\alpha)$
of all $\alpha$-twisted coherent sheaves, where $\alpha$ is a
fixed Brauer class on the surface $X$.

The complexity of a  K3 category $\cat{T}$ is reflected by its
group of $k$-linear, exact auto\-equi\-valences $\Aut(\cat{T})$.
It is in fact a highly non-trivial matter to produce non-trivial
autoequivalences and only very few general techniques are known.
Most importantly, there are the so-called spherical twists which
have been studied in detail by Seidel and Thomas in \cite{ST}.

Another way to gain more insight into the structure of a
triangulated category has more recently introduced by Bridgeland.
In \cite{B2} he defines the notion of a stability condition on a
$k$-linear triangulated category $\cat{T}$ which allows one to
decompose any object in terms of distinguished triangles in a
unique way into semistable objects. Very roughly, a stability
condition is a bounded  $t$-structure together with a numerical
function on its heart.

The truly wonderful aspect of Bridgeland's theory, which has been
strongly influenced by work of Douglas on $\Pi$-stability in the
context of mirror symmetry, is that the space $\Stab(\cat{T})$ of
all stability conditions on a triangulated category $\cat{T}$
admits a natural topology. In general, however, not much is known
about the structure of $\Stab(\cat{T})$. In particular, basic
questions like whether the space is non-empty, connected or
simply-connected are very hard to come by.

In the geometric situation, i.e.\ for the bounded derived
category $\Db(X)$ of all coherent sheaves on a complex variety
$X$, the space $\Stab(X):=\Stab(\Db(X))$ of all stability
conditions is the algebraic version of the usual K\"ahler or
ample cone. Very few stability conditions, i.e.\ points in
$\Stab(X)$, have explicitly been described so far. A complete
description of $\Stab(X)$ for a complex projective (or compact)
manifold $X$ has only been achieved when $X$ is a  curve
\cite{B2,Ma,O}. For partial results for higher-dimensional Fano
manifolds see \cite{A,Ma}.

For compact Calabi--Yau manifolds of dimension at least three not
much is known, not even whether stability conditions always exist.
(Non-compact Calabi--Yau manifolds have been dealt with e.g.\ in
\cite{B4}.) In dimension two, i.e.\ for K3 and abelian surfaces,
Bridgeland constructs in \cite{B} explicit examples of stability
conditions on $\Db(X)$ and investigates in detail one connected
component $\Stab^\dag(X)$ of $\Stab(X)$. Maybe the most fascinating
aspect of the  intriguing paper \cite{B} is a conjecture that
describes $\Stab(X)$ (or rather $\Stab^\dag(X)$) as a universal cover
of an explicit period domain such that the subgroup of
$\Aut(\Db(X))$ of cohomologically trivial autoequivalences acts
naturally as the group of deck transformations. This paper proves
a stronger version of this conjecture for generic twisted K3
surfaces (and abelian surfaces) and generic non-projective K3
surfaces. This provides further evidence for Bridgeland's
conjecture.

The group of autoequivalences of a triangulated category $\cat{T}$
naturally acts on the space of stability conditions
$\Stab(\cat{T})$ and usually the complexity of the group is
reflected by the topology of the quotient. For the special case of
K3 categories it seems that solely spherical twists are
responsible for the rich structure of the group of
autoequivalences, e.g.\ for braid group actions as described in
\cite{ST}. Studied from the point of view of stability conditions,
spherical twists lead to a complicated topological structure of
the quotient of the space of stability conditions by the action of
the group of autoequivalences. This paper confirms this belief by
showing that in the absence of spherical objects, the group of
autoequivalences and the space of stability conditions become
manageable. As generic twisted projective K3 surfaces give rise
to K3 categories not containing any spherical objects, we call
these categories generic K3 categories.

The next more complicated case would be a K3 category with only
one (up to shift) spherical object. In this case one expects that
the group of autoequivalences is generated by the associated
spherical twist, the usual shift functor and a few other obvious
autoequivalences. Geometrically this situation is realized by the
bounded derived category $\Db(X)$ of a generic non-projective K3
surface. Indeed, any line bundle on $X$ defines a spherical
object in $\Db(X)$ and if $X$ is generic non-projective, then the
trivial line bundle is the only one there is.

\medskip

Here are the main results proved in the present paper.

\smallskip

\noindent {\bf Theorem 1.} \emph{Let $X$ be a complex projective
K3 or abelian surface with a Brauer class $\alpha\in
H^2(X,\ko_X^*)$ such that the abelian category of
$\alpha$-twisted coherent sheaves $\coh(X,\alpha)$ does not
contain any spherical objects. Then the space $\Stab(X,\alpha)$
of all locally finite, numerical stability conditions on
$\Db(X,\alpha):=\Db(\coh(X,\alpha))$ admits only one connected
component of maximal dimension. Moreover, this component is
simply-connected. (Theorem \ref{thm:generictwisted})}

\smallskip

It can be shown that in the period domain of all twisted K3
surfaces $(X,\alpha)$ those satisfying the assumption of Theorem 1
form a dense subset (Corollary \ref{cor:nospher}). We expect
$\Stab(\Db(X))$ to be connected, but a technical problem prevents
us, for the time being, from stating it in this generality. See
Remark \ref{rem:excuse} for comments.

\medskip

\noindent {\bf Theorem 2.} \emph{Suppose $(X,\alpha)$ satisfies
the assumption of Theorem 1. Then a Hodge isometry $\widetilde
H(X,\alpha,\ZZ)\cong\widetilde H(Y,\beta,\ZZ)$ can be lifted to an
equivalence $\Db(X,\alpha)\cong\Db(Y,\beta)$ if and only if it
preserves the orientation of the four positive directions.}

\emph{Moreover, there exists a natural surjection
$\Aut(\Db(X,\alpha))\twoheadrightarrow\Aut^+(\widetilde
H(X,\alpha,\ZZ))$ onto the group of orientation preserving Hodge
isometries. The kernel is spanned by the double shift $E\mapsto
E[2]$. (Theorem \ref{thm:genAut}, Corollary \ref{cor:twistder})}

\medskip

The group of autoequivalences is certainly more complicated for
an arbitrary (twisted) K3 surface, mainly due to the existence of
various spherical objects. The first assertion, however, is
expected to hold always, but even for untwisted K3 surface it is
still open whether the Hodge isometry induced by an equivalence
preserves the natural orientation of the positive directions.

\medskip

\noindent {\bf Theorem 3.}  \emph{Suppose $X$ is a complex K3
surface with $\Pic(X)=0$. Then the space $\Stab(X)$ of all
locally finite, numerical stability conditions on
$\Db(X):=\Db(\coh(X))$ is connected and simply-connected.
(Theorem \ref{thm:classgen})}

\medskip

Clearly, a K3 surface satisfying the hypothesis of the theorem
cannot be projective. The abelian category $\coh(X)$ and its
bounded derived category $\Db(X)$ are, in an imprecise sense,
smaller than in the projective situation. However, although there
are no non-trivial line bundles, one still has the point sheaves
$k(x)$ and many vector bundles (necessarily with trivial
determinant). In fact, stable bundles with trivial determinant on
projective K3 surfaces survive the deformation to the generic K3
surface. So, Theorem 3 really describes $\Stab(X)$ in a geometric
highly non-trivial case.

The techniques used to prove Theorem 1 and 3 can be combined to
cover also the case of non-projective twisted K3 surfaces which
are generic in an appropriate sense.

\medskip

\noindent{\bf Theorem 4.} \emph{For a K3 surface as in Theorem 3
one has $\Aut(\Db(X))=\ZZ\oplus\ZZ\oplus\Aut(X)$. (Proposition
\ref{prop:autgennonproj})}

\medskip

Here, $\Aut(\Db(X))$ denotes the group of all autoequivalences of
Fourier--Mukai type, which in general, and this is a difference to
the projective case, is smaller than the full group of all
autoequivalences. The two free factors are spanned by the
spherical twist induced by the trivial line bundle respectively
the shift functor.

\medskip

 We believe that these results for generic twisted and
generic non-projective untwisted K3 surfaces will eventually lead
to a better understanding of stability conditions and
autoequivalences for $\Db(X)$ when $X$ is an arbitrary
(projective) K3 surface.
%For the moment there are a number of
%technical issues related to the deformation theory of categories,
%but we hope to come back to this in the future.

\medskip

The plan of the paper is as follows. In Section \ref{sect:GenCat}
we discuss abstract K3 categories and the role of rigid,
semi-rigid, and spherical objects. In particular, it is shown that
in the absence of spherical objects rigid objects do not exist
(Proposition \ref{prop:stablefact}) and semi-rigid ones are
automatically stable with respect to any stability condition
(Corollary \ref{lem:skyscap}). If there is only one spherical
object for which the spherical twist can be defined, then any
semi-rigid object is stable up to the action of the corresponding
spherical twist (Proposition \ref{prop:onlyone}). The discussion
is later used in the geometric case, but should be applicable in
other situations as well.

Section \ref{sect:twist} is devoted to twisted projective K3 and
abelian surfaces. We first outline how the arguments of \cite{B}
can be adapted to the twisted case (Section \ref{subsect:twist}).
Although dealing at the same time with B-field lifts of Brauer
classes and B-fields complexifying the K\"ahler cone makes the
situation more technical,  most of Bridgeland's arguments go
through unchanged and we will therefore be brief. In Section
\ref{subsect:twistnosph} we apply the general results of Section
\ref{sect:GenCat} to the case of a projective K3 or abelian
surface endowed with a Brauer class $\alpha$ such that
$\Db(X,\alpha)$ does not contain any spherical object. The main
result is Theorem 1 above. Section \ref{subsect:twisauto} presents
two approaches to a complete description of $\Aut(\Db(X,\alpha))$
in the generic case (Theorem 2). The first one follows
Bridgeland's, i.e.\ studying $\Aut(\Db(X,\alpha))$ via its action
on $\Stab(X,\alpha)$, whereas the second one works more directly
by showing that up to shift any Fourier--Mukai kernel is
isomorphic to twisted sheaves.

Section \ref{sect:nonalgK3} deals with K3 surfaces $X$ which do
not admit any non-trivial line bundles. Although these surfaces
are certainly not projective, most aspects of the general theory
go through. As we prove in Section \ref{sect:gennonpro}, the only
spherical object in the derived category $\Db(X)$ is up to shift
the trivial line bundle and thus the discussion of Section
\ref{sect:GenCat} applies. By an ad hoc argument we prove in
Section \ref{sect:Class} that $\Stab(X)$ is connected and
simply-connected (Theorem 3). A more detailed description of
$\Stab(X)$ as the universal covering  of an appropriate period
domain, whose definition differs from the one in the projective
case, has been given in \cite{M}. In order to determine the group
of autoequivalences one could, as for generic twisted K3 surfaces,
imitate Bridgeland's approach and study the action of
$\Aut(\Db(X))$ on $\Stab(X)$ or describe all Fourier--Mukai
kernels explicitly. In Section \ref{subsect:autononproj} we do the
latter and prove Theorem 4.

\smallskip

We would like to point out that similar questions are dealt with in
recent articles by Okada \cite{Ok2} and Ishii, Ueda, and Uehara
\cite{IUU}.

\medskip

\noindent{\bf Acknowledgements.} We would like to thank T.\
Bridgeland for patiently answering our questions and S.\
Meinhardt for discussions related to the results presented in
Section \ref{sect:Class}. The project started when P.S. was
visiting the Max-Planck Institute, whose hospitality is gratefully
acknowledged.

%%%%%%%%%%%%%%%%%%%%%%%%%%%%%%%%%%%%%%%%%%%%%%%%%%%%%%%%%%%%%%%%%%%%%%%%%
\section{K3 Categories: Spherical and (semi-)rigid
objects}\label{sect:GenCat}

This section contains all the abstract arguments of the paper that
are valid for all triangulated categories for which a Serre
functor is provided by the double shift. Although we call these
categories K3 categories, K3 surfaces actually never occur. We
recall the definition of a stability condition and study the
space of all stability conditions on a K3 category without (or
just a few) spherical objects.

\bigskip
%%%%%%%%%%%%%%%%%%%%%%%%%%%%%%%%%%%%%%%%%%%%%%%%%%%%%%%%%%%%%%%%%

\subsection{Stability conditions, phases}
~
\smallskip

\begin{definition}
A \emph{K3 category} is a triangulated category $\cat{T}$
such that the double shift $E\mapsto E[2]$ is a Serre functor.
\end{definition}

All our categories are tacitly assumed to be linear over some
base field $k$ of characteristic $\ne2$ (usually $k=\CC$) and of
finite type, i.e.\ for any two objects $E,F\in \cat{T}$ the
$k$-vector space $\Hom^*(E,F)=\bigoplus \Hom(E,F[i])$ is
finite-dimensional.

We shall use the following short hands: $$(E,F)^i:=\dim
\Hom(E,F[i]) ~{\rm and}~ (E,F)^{\leq i}:=\sum_{j\leq
i}(-1)^j(E,F)^j.$$ Similarly, $(E,F)^{< i}:=(E,F)^{\leq i-1}.$
Thus, if $\cat{T}$ is a K3 category, then $(E,F)^i=(F,E)^{2-i}$.
Moreover, $(E,E)^1$ is always even, as Serre duality defines
a symplectic pairing on it (see \cite{vdB}).

Recall that for two objects $E,F\in \cat{T}$ one defines
$$\chi(E,F):=\sum (-1)^i(E,F)^i.$$
Clearly, $\chi(E[i],F)=\chi(E,F[-i])=(-1)^i\chi(E,F)$ and if
$\cat{T}$ is a K3 category, then $\chi(E,F)$ is symmetric.

Two objects $E_1,E_2\in \cat{T}$ are called \emph{numerically
equivalent}, $E_1\sim E_2$, if $\chi(E_1,F)=\chi(E_2,F)$ for all
$F\in \cat{T}$. In particular, $E$ is numerically trivial if
$E\sim 0$. Note that for all $E\in \cat{T}$ one has $E\sim E[2]$.
The \emph{numerical Grothendieck group} of $\cat{T}$ is defined as
${\cal N}(\cat{T}):=K(\cat{T})/\sim$.

\begin{definition}\label{def:stabcond} ({\bf Bridgeland})
A \emph{stability condition} $\sigma=(Z,\kp)$ on a $k$-linear
triangulated category $\cat{T}$ consists of a linear map
$$\xymatrix{Z:\kn(\cat{T})\ar[r]&\CC}$$ (the \emph{central charge}) and full
additive subcategories $$\kp(\phi)\subset\cat{T}$$ for each
$\phi\in\RR$ satisfying the following conditions:\\
a) If $0\ne E\in\kp(\phi)$, then $Z(E)=m(E)\exp(i\pi\phi)$ for
some
$m(E)\in\RR_{>0}$.\\
b) $\kp(\phi+1)=\kp(\phi)[1]$ for all $\phi$.\\
c) $\Hom(E_1,E_2)=0$ for all $E_i\in\kp(\phi_i)$ with
$\phi_1>\phi_2$.\\
d) Any $0\ne E\in\cat{T}$ admits a \emph{Harder--Narasimhan
filtration} given by a collection of distinguished triangles
$E_{i-1}\to E_i\to A_i$ with $E_0=0$ and $E_n=E$ such that
$A_i\in\kp(\phi_i)$ with $\phi_1>\ldots>\phi_n$.
\end{definition}

This is what Bridgeland calls a numerical stability condition,
but we will only work with stability conditions of this type.

The category $\kp(\phi)$ is the category of \emph{semistable}
objects of phase $\phi$. The objects $A_i$ in condition d) are
called the \emph{semistable factors} or \emph{Harder--Narasimhan
factors} of $E$. They are unique up to isomorphism. The minimal
objects of $\kp(\phi)$ are called \emph{stable} of phase $\phi$.
Note that a non-trivial homomorphism between stable objects of
the same phase is an isomorphism. In particular, for any stable
$E$ the endomorphism algebra $\End(E)$ is a division algebra.
Also, if $E$ is stable, then $(E,E)^i=0$ for $i<0$ and $i>2$.

Bridgeland shows in \cite[Prop.\ 5.3] {B2} that a stability
condition can equivalently be described by a bounded $t$-structure
and a centered stability function $Z:\ka\to\CC$ on its heart $\ka$
which has the Harder--Narasimhan property. More precisely, one has
$\ka=\kp((0,1])$, i.e.\ $\ka$ contains all objects $E\in\cat{T}$
whose semistable factors $A_i$ satisfy $\phi(A_i)\in(0,1]$. We
call $\ka$ also the heart of the stability condition $\sigma$.

Following Bridgeland, all stability conditions we will consider
are locally finite \cite{B}. In particular, any semistable object
$E\in\kp(\phi)$ admits a finite Jordan--H\"older filtration, i.e.\
a finite filtration $E_0\subset\ldots\subset E_n=E$ with stable
quotients $E_{i+1}/E_i\in\kp(\phi)$.

The space of all locally finite
stability conditions is denoted $\Stab(\cat{T})$.
For later use we state the following observation, which will be
needed to ensure equality of phases of a distinguished family of
stable objects.

\begin{lem}\label{prop:phase}
Suppose $\cat{T}$ is a K3 category endowed with a stability condition
$\sigma$ and suppose $E_1\sim E_2$ are numerically equivalent $\sigma$-stable
objects. If there exists a $\sigma$-stable object
 $E_0\in \cat{T}$ with $$(E_0,E_1)^0\ne0\ne(E_0,E_2)^0,$$ then
$$\text{either}~~\phi(E_1)=\phi(E_2)~~\text{or}~~
E_1\cong E_2[\pm2].$$
\end{lem}

\begin{proof}
Let $\phi_i:=\phi(E_i)$. By assumption and Serre duality there
exist non-trivial homomorphisms $E_0\to E_i\to E_0[2]$, which due
to stability implies $\phi_0\leq\phi_i\leq\phi_0+2$.

As $E_1\sim E_2$, the values of the stability function satisfy
$Z(E_1)=Z(E_2)$ and hence $\phi_1=\phi_2+2\ell$ for some integer
$\ell$. Thus, either $\phi_1=\phi_2$ or one has one of the
following possibilities i) $\phi_1=\phi_0=\phi_2-2$ or ii)
$\phi_2=\phi_0=\phi_1-2$. Due to the stability of all three
objects, the latter two cases lead to i) $E_1\cong E_0\cong
E_2[-2]$ respectively ii)  $E_2\cong E_0\cong E_1[-2]$.
\end{proof}

\begin{remark}\label{rem:avoid} In the geometric context, numerical equivalence of two
objects (e.g.\ sheaves) $E_1$ and $E_2$  is often caused by the
existence of a  `flat family of objects (sheaves)' connecting
$E_1$ and $E_2$. In this context, for $E_1$ and $E_2$ stable objects, $\phi(E_1)=\phi(E_2)$ follows
from the `openness of stability' as proved in \cite[Prop.\
3.3.2]{AP}. So in fact the above proposition could be avoided
altogether in the geometric situation we shall be interested in
later.
\end{remark}

%\begin{cor}\label{cor:phase}
%Let ${\mathcal A}$ be a $k$-linear
%abelian category such that $\cat{T}:=\Db({\mathcal A})$
%is a K3 category. Suppose $E_i\in{\mathcal A}$, $i=0,1,2$,
%are objects with
%$E_1\sim E_2$ and $\chi(E_0,E_1)\ne0$ (or equivalently $\chi(E_0,E_2)\ne0$).
%
%If for a stability condition $\sigma$ all three objects $E_i$ are
%$\sigma$-stable, then $$\phi(E_1)=\phi(E_2).$$
%\end{cor}

%\begin{proof}
%Suppose $\chi(E_0,E_1)=\chi(E_0,E_2)>0$. Then there are, up to
%changing the order of $E_1$ and $E_2$, the following two cases: i)
%$(E_0,E_1)^0\ne0$ and $(E_0,E_2)^0\ne0$ or ii) $(E_0,E_1)^0\ne0$
%and $(E_0,E_2)^2\ne0$. In case i), the proposition applies
%directly, for $E_1\not\cong E_2[\pm2]$ and yields
%$\phi(E_1)=\phi(E_2)$. In case ii), the proposition yields
%$\phi(E_1)=\phi(E_2[2])=\phi(E_2)+2$ or $E_1\cong E_2[2][\pm2]$.
%The latter is only possible if $E_1\cong E_2$ and then of course
%$\phi(E_1)=\phi(E_2)$. If, however, $\phi(E_1)=\phi(E_2)+2$, one
%obtains a contradiction from $\phi(E_2)\leq\phi(E_0)\leq
%\phi(E_1)+2$???????

%If $\chi(E_0,E_1)=\chi(E_0,E_1)<0$, then $(E_0,E_i)^1\ne$ for
%$i=1,2$ and the proposition applies directly.
% ~~{\rm or}~~
%\chi(E_0,E_1)=\chi(E_0,E_2)<0.$$ In the first case,
%$(E_0,E_1)\ne0\ne(E_0,E_2)$ and one can apply the proposition. In
%the second case $(E_0[-1],E_1)\ne0\ne (E_0[-1],E_2)$ and the
%proposition applies again with $E_0$ replaced by $E_0[-1]$.
%\end{proof}
%%%%%%%%%%%%%%%%%%%%%%%%%%%%%%%%%%%%%%%%%%%%%%%%%%%%%%%

\bigskip

\subsection{Spherical, rigid, and semi-rigid
objects}\label{subsect:sphercial} ~
\smallskip

Following the standard terminology in \cite{Mu} and \cite{ST}, we
give the following definition:

\begin{definition}\label{def:sohrig}
i) An object $E\in\cat{T}$
is \emph{rigid} if $(E,E)^1=0$;

ii) An object $E\in\cat{T}$
is \emph{semi-rigid} if $(E,E)^1=2$;

iii)  An object $E\in\cat{T}$ is \emph{spherical} if $(E,E)^i=1$
for $i=0,2$ and otherwise zero.

iv) An object $E\in\cat{T}$ is \emph{quasi-spherical} if
$(E,E)^i=0$ for $i\ne0,2$ and $\End(E)$ is a division algebra.
\end{definition}

\begin{remark}
Since any finite-dimensional division algebra over an
algebraically closed field $k$ is isomorphic to $k$ itself, the
notions of quasi-spherical and spherical coincide in the case of
a $k$-linear category over an algebraically closed field $k$.

In fact, for the results mentioned in the introduction it would be
enough to work over an algebraically closed field, but for the
sequel to this paper it is important to consider the more general
situation as well. Indeed, in \cite{HMS2} we  use the techniques
in the present paper to prove that any autoequivalence of the
derived category of a smooth projective K3 surface induces an
orientation preserving Hodge isometry on cohomology. Our argument
involves the passage to a triangulated category which is linear
over a field which is not algebraically closed.\end{remark}

In order to control how these notions behave in distinguished
triangles, the following lemma will be crucial. It is the
analogue of Lemma 5.2 in \cite{B}, but the main idea goes back to
Mukai's paper \cite{Mu}.

\begin{lem}\label{lem:Mukaigen}
 Consider in $\cat{T}$ the distinguished triangle
$$\xymatrix{
A\ar[r]^i&E\ar[r]^j&B\ar[r]^\delta&A[1]}$$ such that
\[
(A,B)^{r}=(B,B)^{s}=0~\text{for~} r\leq0\text{~and~} s<0.
\]
Then
\[
(A,A)^1+(B,B)^1\leq(E,E)^1.
\]
\end{lem}

\begin{proof} First of all we note that
$\chi(E,E)=\chi(A,A)+\chi(B,B)+2\chi(A,B)$. Due to our hypotheses
and using Serre duality this can be rewritten as
\begin{eqnarray}\label{eqn:1}
(A,A)^1+(B,B)^1=(E,E)^1+2\underbrace{\left((A,A)^{\leq
0}+(B,B)^0-(E,E)^{\leq 0}+(B,A)^{\leq 1}\right)}_{=N}.
\end{eqnarray}
The result is proved once we show that $N\leq0$.

Since $\Hom(A,B)=\Hom(A,B[-1])=0$, for any $f\in\Hom(E,E)$, there
exist unique $g\in\Hom(A,A)$ and $h\in\Hom(B,B)$ making the
following diagram commutative:
\begin{eqnarray}\label{eqn:2}
\xymatrix{A\ar[r]^{i}\ar[d]^{g}&E
\ar[r]^{j}\ar[d]^{f}&B\ar[r]^-{\delta}\ar[d]^{h}&A[1]\ar[d]^{g[1]}\\
A\ar[r]^{i}& E\ar[r]^{j}&B\ar[r]^-{\delta}&A[1].}
\end{eqnarray}
In particular, we get a map
$\Hom(E,E)\to\Hom(A,A)\oplus\Hom(B,B)$.

If $g=h=0$, then $f$ factorizes through a morphism $f_1:E\to A$.
Moreover, since $i\circ f_1\circ i=i\circ g=0$, the morphism
$f_1\circ i$ can be factorized further through some $f_2:A\to
B[-1]$ which must be trivial due to our assumption.
 Thus $f_1\circ i=0$ and hence $f_1$ admits a
factorization through $f_3:B\to A$. This discussion leads to an
exact sequence
$$\xymatrix{
\Hom(B,A)\ar[r]&\Hom(E,E)\ar[r]^-{\tau}&\Hom(A,A)\oplus\Hom(B,B).}$$

Consider now the morphism
$\eta:\Hom(A,A)\oplus\Hom(B,B)\to\Hom(B,A[1])$ which we define as
$\eta((g,h))=\delta\circ h-g[1]\circ\delta$. Suppose that
$(g,h)=\tau(f)$, for some $f\in\Hom(E,E)$. Then $\eta((g,h))=0$
and it is very easy to check that the previous exact sequence can
be completed to the following:
$$\xymatrix{
\Hom(B,A)\ar[r]^-\zeta&\Hom(E,E)\ar[r]^-\tau&\Hom(A,A)\oplus\Hom(B,B)\ar[r]^-\eta&
\Hom(B,A[1]).}$$

Next, take $g_1:B\to A$ such that $\zeta(g_1)=0$. By definition
this means that $i\circ g_1\circ j=0$ and then $g_1\circ j$
factorizes through some $g_2:E\to B[-1]$, which must be trivial.
Indeed, otherwise $g_2$ would  give rise to a non-trivial $B\to
B[-1]$ or a non-trivial $A\to B[-1]$, the existence of which is in
both cases excluded by assumption. Hence, $g_1\circ j=0$ and,
therefore, $g_1$ factorizes through some $g_3:A[1]\to A$. Thus we
obtain the exact sequence
$$\xymatrix@=14pt{\Hom(A[1],A)\ar[r]^-\theta&\Hom(B,A)\ar[r]^-\zeta&
\Hom(E,E)\ar[r]^-\tau&\Hom(A,A)\oplus\Hom(B,B)\ar[r]^-\eta&\Hom(B,A[1]).}$$

Reasoning as before, consider $g'\in\Hom(A[1],A)$ such that
$\theta(g')=0$. Then $g'\circ\delta=0$ and there exists
$g'_1:E[1]\to A$ such that $g'=g'_1\circ i[1]$. Using the exact
sequence
$$\xymatrix{
\Hom(E[1],B[-1])\ar[r]&\Hom(E[1],A)\ar[r]&\Hom(E[1],E)\ar[r]&
\Hom(E[1],B)}$$ and the assumption $(E[1],B)^m=0$ for $m\leq0$ we
see that $g'_1$ induces a unique $g'_2:E[1]\to E$. Thus the
sequence
$$\xymatrix{
\Hom(E[1],E)\ar[r]&\Hom(A[1],A)\ar[r]^\theta&\Hom(B,A) }$$ is
exact.

Iterating the previous argument we obtain the exact sequence
\begin{eqnarray}\label{eqn:exseq}
\begin{array}{c}
\cdots\mor\Hom^{-m}(E,E)\mor\Hom^{-m}(A,A)\mor\Hom^{-m+1}(B,A)\mor\cdots\\
\cdots\mor\Hom^{0}(E,E)\mor\Hom^{0}(A,A)\oplus\Hom^{0}(B,B)\mor[\eta]
\Hom^{1}(B,A).
\end{array}
\end{eqnarray}
Therefore, $-(E,E)^{\leq 0}+(A,A)^{\leq 0}+(B,B)^{0} + (B,
A)^{\leq 0} \leq (B,A)^1$, i.e.\ $N\leq 0$.
\end{proof}

\begin{remark}\label{rem:cond}
If the base field is not algebraically closed, we have to
introduce the following technical condition which ensures that the
theory goes through smoothly. Situations where this condition does
not hold seem rather pathological. In any case, in all
applications dealt with in this paper and its sequel, it can be
easily verified. See e.g.\ Proposition \ref{rem:mod}.

\smallskip

\begin{itemize}
\item[(${\bf \ast}$)] {\it The category $\cat{T}$ does not contain any object $B\in\cat{T}$ such that $\End(B)$ is a division algebra,  $(B,B)^0=(B,B)^1$ and $(B,B)^i=0$ for any integer $i<0$.}
\end{itemize}

\smallskip

Due to the existence of the symplectic structure on
$\Hom(B,B[1])$, we know that it is an even-dimensional $k$-vector
space. So, if $k$ is algebraically closed, this condition is
automatic, but in general  $\Hom(B,B[1])$ need not be
even-dimensional over $\End(B)$.
\end{remark}

\begin{prop}\label{prop:stablefact} Suppose $\cat{T}$ is a K3 category
satisfying condition {\rm ($\ast$)} and let  $\ss\in\Stab(\cat{T})$ be a
stability condition.

{\rm i)} If $E\in\cat{T}$ is rigid, then all $\ss$-stable factors
of $E$ are quasi-spherical (and thus spherical if $k$ is
algebraically closed).

{\rm ii)} If $E\in\cat{T}$ is semi-rigid, then all $\ss$-stable
factors of $E$ are either quasi-spherical or semi-rigid. In fact,
there exists at most one semi-rigid stable factor, which moreover
occurs with multiplicity one.
\end{prop}

\begin{proof} We start out with the following technical
observation, which is irrelevant if $k$ is algebraically closed.
Suppose that the endomorphism algebra $\End(B)$ of an object
$B\in\cat{T}$ is a division algebra. Then $\Hom(B,B[1])$ is a
vector space over $\End(B)$. Note however that multiplication
from the right and from the left do not necessarily commute. In
particular, as an abstract $\End(B)$-vector space one has
$\Hom(B,B[1])\cong\End(B)^{\oplus r}$ for some $r$ and, therefore,
$(B,B)^1=r(B,B)^0$.
%Furthermore, using the symplectic structure
%$\Hom(B,B[1])\times\Hom(B,B[1])\to\End(B)^*$, one proves
%$r\equiv0(2)$.

Let $A_1,\ldots,A_n$ be the semistable factors of an object
$E\in\cat{T}$. Then Lemma \ref{lem:Mukaigen} implies that
$\sum(A_i,A_i)^1\leq (E,E)^1$.

i) Suppose now that $E$ is rigid. Then the $A_i$ are as well
rigid. Thus, since an object that is stable and rigid is
quasi-spherical, it suffices to show that the stable factors of a
semistable rigid $A\in\cat{T}$ are also  rigid.

Let us first assume that the stable factors of $A$ are all
isomorphic to the same object $B$.  Then
$\chi(A,A)=\ell^2\chi(B,B)$ for some positive integer $\ell$. By
rigidity of $A$ one has $\chi(A,A)=2(A,A)^0>0$ and for $B$ one has
$\chi(B,B)=2(B,B)^0-(B,B)^1$. Since $B$ is stable, $\End(B)$ is a
division algebra and the above observation applies. Hence,
$\chi(B,B)=(2-r)(B,B)^0$ for some $r\geq0$ and
$(2-r)\ell^2=2(A,A)^0/(B,B)^0>0$ then shows that $r=0$ or $r=1$,
but the latter is excluded by condition ($\ast$).

Suppose now that $A$ has at least two non-isomorphic stable
factors. Then there exists a distinguished triangle
\begin{equation}
\label{JH}\xymatrix{ C\ar[r]& A\ar[r]& D\ar[r]&C[1]}
\end{equation}
 with $C$ and $D$
semistable of the same phase and such that $(C,D)^i=0$, for any
$i\leq 0$, and such that all stable factors of $D$ are isomorphic
to each other. Again by Lemma \ref{lem:Mukaigen},  $C$ and $D$ are
rigid. Then by the above argument the stable factors of $D$ are
quasi-spherical  and one continues with $C$.

ii) If $E$ is semi-rigid, then all its semistable factors will be
either rigid or semi-rigid. Due to i), the stable factors of a rigid
semistable $A$ will all be quasi-spherical. If $A$ is semi-rigid and
semistable with just one type of stable factor $B$, then as before,
$\chi(A,A)=\ell^2\chi(B,B)$. This time one has only the weak
inequality $\chi(A,A)\geq0$, which yields $(2-r)(B,B)^0\geq0$, i.e.\
$r=0$, $r=1$, or $r=2$. If $r=2$, then $\chi(B,B)=0$ and hence
$\chi(A,A)=0$. As $A$ is semi-rigid, the latter implies that $A$ is
simple and hence $A=B$. Again the case $r=1$ is excluded by
condition ($\ast$). If $r=0$, then $B$ is rigid and thus $A$ would
be direct sum of copies of $B$ which is not the case.

If  $A$ has at least two non-isomorphic stable factors, then use
the distinguished triangle (\ref{JH}) as above. At most one of the
two objects $C$ or $D$ can be semi-rigid. More precisely, $A$ has
only one semi-rigid stable factor which can only occur in either
$C$ or $D$. Combined with the discussion before one sees that all
stable factors of $E$ are spherical up to at most one which is
semi-rigid and occurs with multiplicity one.
\end{proof}

For $\sigma=(Z,\kp)\in\Stab(\cat{T})$ and a $\sigma$-semistable $0\neq E\in\kp(1)$, an analogous result is proved in \cite[Lemma 12.2]{B}. In the applications the following corollary will be combined with
Lemma \ref{prop:phase}.

\begin{cor}\label{lem:skyscap}
Let $\cat{T}$ be a K3 category satisfying {\rm ($\ast$)} and not
containing any quasi-spherical objects. If $E\in\cat{T}$ is
semi-rigid, then $E$ is stable with respect to any stability
condition.\qqed\end{cor}

The above discussion can be generalized to yield  a relation
between $(E,E)^1$ and the number of (semi)stable factors of $E$.
For this we denote by $\ell_{\rm HN}(E)$ the number of semistable
(Harder--Narasimhan) factors of an object $E$ and, similarly, by
$\ell_{\rm JH}(E)$ the number of stable (Jordan--H\"older) factors
counted without(!) multiplicity. Then one has

\begin{cor}
Let $\cat{T}$ be a K3 category satisfying {\rm ($\ast$)} and not
containing any quasi-spherical objects. Then $\ell_{\rm
HN}(E)\leq\ell_{\rm JH}(E)\leq(1/2)(E,E)^1$.\qqed
\end{cor}

\bigskip

%%%%%%%%%%%%%%%%%%%%%%%%%%%%%%%%%%%%%%%%%%%%%%%%%%%%%%%%%%%%%%%%%
\subsection{Only one spherical object}\label{sect:onlyone}
~

\smallskip

So far we have dealt with K3 categories not containing any
spherical objects. Let us now pass to the slightly more
complicated situation when there exists a spherical object $E\in
\cat{T}$ which is unique up to shift.

Let us first explain how to construct inductively indecomposable
rigid objects  starting with an arbitrary spherical object $E$.
This goes as follows: The unique morphism
$$\xymatrix{E\ar[r]&E[2]}$$
can be completed to a distinguished triangle
$$\xymatrix{E_2\ar[r]&E\ar[r]&E[2]\ar[r]&E_2[1].}$$
One checks easily that $(E,E_2)^2=1$ and that therefore there
exists a natural distinguished triangle
$$\xymatrix{E_3\ar[r]&E\ar[r]&E_2[2]\ar[r]&E_3[1].}$$
Due to the following lemma this procedure can be iterated and
produces distinguished triangles
\begin{equation}\label{Indu}
\xymatrix{E_{n+1}\ar[r]&E\ar[r]&E_n[2]\ar[r]&E_{n+1}[1].}
\end{equation}
Thus we obtain a sequence of rigid objects
$E_1:=E,E_2,E_3,\ldots$.

\begin{lem}\label{lem:calc}
If $m \leq n$ are two positive integers, then
\begin{equation*}
(E_m,E_n)^i=
              \begin{cases}
                          1 & \text{if} \ i= -n+1, \ldots , -n+m \\
                          1 & \text{if} \ i= 2, \ldots , m+1 \\
                          0 & \text{otherwise}.
              \end{cases}
\end{equation*}\end{lem}

\begin{proof}
Applying $\Hom(E,-)$ to (\ref{Indu}) we obtain the long exact
sequence
$$\xymatrix{\ldots\ar[r]&\Hom^i(E,E_n[1])\ar[r]&\Hom^i(E,E_{n+1})\ar[r]
&\Hom^i(E,E)\ar[r]&\ldots}$$ and so by induction
\begin{equation*}
(E,E_n)^i=
              \begin{cases}
                          1 & \text{if} \ i= -n+1 \\
                          1 & \text{if} \ i=2 \\
                          0 & \text{otherwise}.
              \end{cases}
\end{equation*}

Applying the exact sequence \eqref{eqn:exseq} in the proof of
Lemma \ref{lem:Mukaigen} to the distinguished triangle
$E_{n-1}[1]\to E_n\to E$, we obtain the exact sequence
$$\xymatrix@=14pt{\ldots\ar[r]&\Hom^{-i}(E,E_{n-1}[1])\ar[r]&
\Hom^{-i}(E_n,E_n)
\ar[r]&\Hom^{-i}(E_{n-1}[1],E_{n-1}[1])\ar[r]&\ldots,}$$
 if $i > 0$ and
\begin{multline*}
\xymatrix@=14pt{\ldots\ar[r]&\Hom^0(E,E_{n-1}[1])\ar[r]&\Hom^0(E_n,E_n)\ar[r]&
\Hom^0(E_{n-1}[1],E_{n-1}[1])\oplus\Hom^0(E,E)}\\
\xymatrix@=14pt{\ar[r]^-\eta&\Hom^1(E,E_{n-1}[1]).}
\end{multline*}
Since $\eta$ is surjective, we get the following
\begin{equation*}
(E_n,E_n)^i=
              \begin{cases}
                          1 & \text{if} \ i= -n+1, \ldots , 0 \\
                          1 & \text{if} \ i=2, \ldots , n+1 \\
                          0 & \text{otherwise}.
              \end{cases}
\end{equation*}
Finally, the lemma follows by applying $\Hom (E_m,-)$ to
\eqref{Indu}.
\end{proof}

An immediate consequence of the above is
\begin{cor}
If $E\in\cat{T}$ is a spherical object, then the objects $E_n$
are rigid.
\end{cor}

\begin{prop}\label{rem:mod}
Suppose our K3 category $\cat{T}$ is the bounded derived category
$\Db(\ka)$ of some abelian category $\ka$ and that $\ka$ contains
a spherical object $E\in\ka$ which is the only indecomposable
rigid object in $\ka$. Then:

{\rm i)} Any indecomposable rigid object $F\in\cat{T}$ is
isomorphic to a shift of some $E_n$,

{\rm ii)}  $E$ is up to shift the only quasi-spherical object also in
$\cat{T}$, and

{\rm iii)} $\cat{T}$ satisfies condition {\rm ($\ast$)} in Remark
\ref{rem:cond}.
\end{prop}

\begin{proof} The assertion is proved by induction
on the length of the complex $F$. Since $F$ is rigid, all its
cohomology objects $H^i(F)$ are rigid (see e.g.\ \cite[Lemma
2.9]{BM} or Lemma \ref{lem:Mukaigen}). Suppose $k$ is maximal with $H^k(F)\ne0$ and for
simplicity we assume $k=0$. Then $H^0(F)=E^{\oplus r}$ for some
$r$ and one has the natural distinguished triangle
$$\xymatrix{F_1\ar[r]&F\ar[r]&E^{\oplus r}\ar[r]&F_1[1]}$$
with $H^i(F_1)=0$ for $i\geq0$. In particular, the assumptions of
Lemma \ref{lem:Mukaigen} are satisfied and, therefore, $F_1$ is
again rigid. By induction $F_1=\bigoplus E_i[n_i]^{\oplus r_i}$.

Therefore, $\Hom(E,E_i[n_i+1])\neq 0$, as otherwise $E_i[n_i]$
would be a direct summand of $F$, which contradicts the assumption
that $F$ is indecomposable. By Lemma \ref{lem:calc} the only
non-trivial morphisms $E\to E_i[n_i+1]$ exist for $n_i+1=2$ and
$n_i+1=-i+1$. In the first case, $E_{i+1}$ would be a direct
summand of $F$ and hence $F=E_{i+1}$, for $F$ is indecomposable.
In the latter case $E_i[-i]$ would be a summand of $F_1$, which
contradicts $H^i(F_1)=0$ for $i\geq0$. (Note
$H^i(E_i[-i])=H^0(E_i)\ne0$.)

Let us now verify condition ($\ast$). Suppose $B\in\Db(\ka)$ is an
object such that $\End(B)$ is a division algebra, $(B,B)^0=(B,B)^1$ and $(B,B)^i=0$, for any integer $i<0$. We shall derive a contradiction as follows.
Consider a non-trivial extension $B\to A\to B\to B[1]$. (In fact,
the isomorphism type of $A$ is unique, but we will not need this.)
The boundary morphisms $\End(B)\to\Hom(B,B[1])$ and
$\Hom(B,B[1])\to \Hom(B,B[2])$ are bijective, because of Serre
duality and linearity under the action of $\End(B)$. This shows
$\Hom(B,B)\cong\Hom(B,A)$ and $\Hom(A,A[i])=0$ for $i\ne0,2$. We
leave it to the reader to show that $A$ is also indecomposable.
Hence, $A\cong E$ (up to shift), which yields the contradiction
$2=\chi(E,E)=\chi(A,A)=4\chi(B,B)=4(B,B)^0>2$.
\end{proof}

\begin{prop}\label{prop:stabonly}
Suppose $\cat{T}$ satisfies {\rm ($\ast$)} and $E\in\cat{T}$ is a
spherical object which is up to shift the only  quasi-spherical
object in $\cat{T}$. Then $E$ is stable with respect to any
stability condition on $\cat{T}$ and the stable factors of $E_n$
are $E[n-1],\ldots,E[1], E$.
\end{prop}

\begin{proof}
Consider the stable factors $A_1,\ldots,A_\ell$ of $E$. By
Proposition \ref{prop:stablefact} they are all quasi-spherical and
hence isomorphic to some $E[i]$. %Hence
%$2=\chi(E,E)=\sum\chi(A_i,A_j)=\ell^2\chi(E,E)$, which yields
%$\ell=1$, i.e.\
In particular, $E$ is stable. The description of the stable
factors of the $E_n$ follows immediately.
\end{proof}

\begin{remark} Indecomposable rigid objects in a K3 category with
only one spherical object can be classified in general. The
precise result is the following: Let $\cat{T}$ be a K3 category
satisfying ($\ast$) and containing a spherical object $E$ which is
up to shift the only  quasi-spherical object in $\cat{T}$. Suppose
there exists a stability condition on $\cat{T}$. Then any
indecomposable rigid object in $\cat{T}$ is isomorphic to a shift
of some $E_n$ which are constructed as above.
\end{remark}

%\begin{proof}
%Suppose $F$ is an indecomposable rigid object. Then all its
%stable factors are spherical and hence isomorphic to shifts of
%$E$ (cf.\ Proposition \ref{prop:stablefact}). Therefore, by
%induction it suffices to show that if $F$ is indecomposable and
%given by a `maximal destabilizing' distinguished triangle
%$$\xymatrix{\bigoplus E_i[n_i]^{\oplus r_i}\ar[r]&F\ar[r]&E^{\oplus r}\ar[r]&
%\bigoplus E_i[n_i+1]^{\oplus r_i},}$$ then $F\cong E_k[n_k]$ for
%some $k$ and $n_k$. Here, as in the proof of Proposition
%\ref{prop:stablefact}, one works first with both,  the
%Harder--Narasimhan and Jordan--H\"older filtration. Thus `maximal
%destabilizing' here means that $\Hom(\bigoplus E_i[n_i]^{\oplus
%r_i},E^{\oplus r})=0$.
%
%Clearly, all $E\to E_i[n_i+1]$ must be non-trivial, as otherwise
%$E_i[n_i]$ would be a direct summand of $F$, which contradicts
%the assumption that $F$ is indecomposable. By Lemma
%\ref{lem:calc} the only non-trivial $E\to E_i[n_i+1]$ exist for
%$n_i+1=2$ and $n_i+1=-i+1$. In the first case, $E_{i+1}$ would be
%a direct summand of $F$. In the latter case the non-trivial
%$E_i[-i]\to F$ would contradict the assumption that $\bigoplus
%E_i[n_i]^{\oplus r_i}\to F$ is maximal destabilizing, as
%$(E_i[-i],E[-i])^0\ne0$ shows $\phi(E_i[-i])<\phi(E)$.
%\end{proof}

In all relevant examples of K3 categories a spherical object
$E\in\cat{T}$ gives rise to a \emph{spherical twist}, a particular
exact autoequivalence of $\cat{T}$ (see \cite{ST} or \cite[Ch.\ 8]
{HFM}). In the following, we shall thus assume that the spherical
twist can be constructed, i.e.\ that there exists an
autoequivalence
$$\xymatrix{T_E:\cat{T}\ar[r]&\cat{T},}$$
such that for any object $F\in\cat{T}$ one has a distinguished
triangle
$$\xymatrix{T_E(F)[-1]\ar[r]
&\Hom(E,F[*])\otimes E\ar[r]&F\ar[r]&T_E(F)&}$$ which is
functorial in $F$. Recall that the inverse can be described by a
distinguished triangle
$$\xymatrix{T_E^{-1}(F)\ar[r]&F\ar[r]&E[2]\otimes \Hom(E,F[*])\ar[r]&T_E^{-1}(F)[1].}$$
Also recall that $T_E^k(E)\cong E[-k]$.

\begin{remark}\label{rem:t2not}
i) We shall be interested in the spherical twist $T:=T_E$ applied
to objects $F\in\cat{T}$ with $\sum (E,F)^i=1$. After shifting $F$
(or $E$) appropriately, this reduces to the case $(E,F)^0=1$ and
$(E,F)^i=0$ for $i\ne0$.

Then for all $k$ one has $(E,T^k(F))^k=1$ and $(E,T^k(F))^i=0$ for $i\ne k$.
Indeed, $$(E,T^k(F))^i=(T^k(E)[k],T^k(F))^i=(E[k],F)^i=(E,F)^{i-k}.$$

Using the distinguished triangles $E\to T^k(F)[k]\to
T^{k+1}(F)[k]\to E[1]$, one also proves by induction
$(F,T^k(F))^0=1$ for $k\geq0$.

ii) Suppose the spherical object $E$, $T^n(F)$, and $T^{k}(F)$ are all $\sigma$-stable for
some stability condition $\sigma$ on $\cat{T}$. If as before $\sum(F,E)^i=1$, then $|n-k|<2$.

Clearly, $E$ is also stable with respect $T^n(\sigma)$, so that
we may assume that $n=0$, $k>0$, and $(E,F)^0=1$.

Then using the non-triviality of $(F,T^k(F))^0$,
$(T^k(F)[k],E)^2$, and $(E,F)^0$ one obtains the following
inequality of phases
$$\phi(F)\leq\phi(T^k(F))\leq\phi(E)-k+2\leq\phi(F)-k+2.$$
Hence $k\leq2$. Moreover, if $k=2$ then $\phi(E)=\phi(F)$ and hence $E\cong F$, which
is absurd.

iii) For later use we rephrase the last observation as follows:
Let $E$ be a spherical object in a K3 category $\cat{T}$
satisfying ($\ast$), which induces a spherical twist $T:=T_E:\cat{
T}\cong\cat{T}$. Suppose $E$ is up to shift the only
quasi-spherical object. If $F\in\cat{T}$ with $\sum (E,F)^i=1$ and
$W_F:=\{\sigma\in\Stab(\cat{T})~|~F~{\rm
is}~\sigma\text{-stable}\}$, then
$$T^nW_F\cap T^kW_F=\emptyset$$ if
$|n-k|\geq2$.
\end{remark}

\begin{prop}\label{prop:onlyone}
Let $\cat{T}$ be a K3 category satisfying {\rm ($\ast$)} with a
spherical object $E$ inducing a spherical twist
$T:=T_E:\cat{T}\cong\cat{T}$. Suppose there is no other
quasi-spherical object up to shift. Then for any semi-rigid object
$F$ with $(F,E)^*:=\sum(F,E)^i=1$ and any stability condition
$\sigma$, one finds an integer $n$ such that $T^n(F)$ is stable.
\end{prop}

\begin{proof} By assumption, we may
 assume that $(E,F)^0=(F,E)^2=1$ and $(E,F)^i=(F,E)^{2-i}=0$ for all
$i\ne0$.

Consider the Harder--Narasimhan factors $A_1,\ldots,A_n$ of $F$
with respect to a given stability condition. Then by Proposition
\ref{prop:stablefact} at most one  $A_i$ is semi-rigid and all
others are rigid.

Suppose first that $A_1$ is rigid. Then, again by Proposition
\ref{prop:stablefact}, all stable factors of  $A_1$ are
quasi-spherical and hence isomorphic to some shift of $E$. Since
$A_1$ itself is semistable, all the shifts are the same, say
$E[k]$. In fact $k=0$, for there exists a non-trivial morphism
$A_1\to F$. On the other hand, the indecomposable factors of
$A_1$ are as well rigid and hence isomorphic to some $E_j[n_j]$,
where the $E_j$ are as above. However, due to Proposition
\ref{prop:stabonly} the $E_j$ are semistable only for $j=1$.
Combining both observations yields $A_1\cong E^{\oplus r}$.

Next apply $\Hom(E,-)$ to the distinguished triangle
$$\xymatrix{A_1\ar[r]&F\ar[r]&F'\ar[r]&A_1[1].}$$
Since $(A_1,F')^{-1}=(A_1,F')^0=0$ and $(E,F)^0=1$, this shows $r=1$. Hence,
the  first factor in the Harder--Narasimhan filtration $A_1\to F$
can be interpreted as the natural map $\Hom(E,F)\otimes E\to F$.
Thus, the cone $F'$ is isomorphic to $T(F)$.

Now one continues with $F'$, whose Harder--Narasimhan filtration
has smaller length than the one of $F$. Indeed, $F'$ is again
semi-rigid and $(F',E)^*=(T(F),E)^*=(F,T^{-1}(E))^*=(F,E)^*=1$,
i.e.\ $F'$ satisfies all the assumptions on $F$.

If the first Harder--Narasimhan factor $A_1$ is only semi-rigid,
then we rather work with the last one $A_n$ which is necessarily
rigid. The same reasoning as above shows that $A_n\cong E[2]$ and
this time the last quotient $F\to E[2]$ of the Harder--Narasimhan
filtration can be interpreted as the natural morphism
$F\to\Hom^2(F,E)\otimes E[2]$. Hence, the kernel is isomorphic to
$T^{-1}(F)$. As above, one proceeds then with $T^{-1}(F)$, which
still meets all the requirements.

Thus, there exists an integer $n$ such that $T^n(F)$ is
semistable. (Of course, working through the Harder--Narasimhan
filtration only either $T$ or $T^{-1}$ will occur.)

 We leave it to the reader to repeat the same arguments
once more for appropriate Jordan--H\"older filtrations. Compare
the arguments in the proof of Proposition \ref{prop:stablefact}.
\end{proof}

In the application we shall need stability of $T^n(F)$ for
several $F$ at a time, but with the same $n$. The assertion that
will be needed is the following:

\begin{cor}\label{cor:Tnstableonly}
If $F_1\sim F_2$ are as above with $(E,F_1)^i=(E,F_2)^i$ and
$(F_1,F_2)^i=0$ for all $i$, then there exists an integer $n$ such
that $T^n(F_1)$ and $T^n(F_2)$ are both stable.
\end{cor}

\begin{proof}
As in the proof of the proposition, we shall assume that
$(E,F_1)^0=(E,F_2)^0=1$.
%Since $E$ is stable with respect to any
%%$T^n(\sigma)$ and, if necessary, interchanging the role of $F_1$ a
%nd $F_2$, that $F_1$ and $T^n(F_2)$ are stable for some $n>0$.
 After interchanging the role of
$F_1$ and $F_2$ if necessary, we may assume that $F_1$ and
$T^n(F_2)$ are stable for some $n>0$ (see Remark \ref{rem:t2not},
ii)). We shall show that in fact $F_2$ is already semistable. The
proof of the stability of $F_2$ is similar and uses appropriate
Jordan--H\"older filtrations (instead of Harder--Narasimhan).

Following the arguments in the above proof, we shall assume that
the unique $E\to F_2$ is the first factor in the
Harder--Narasimhan filtration.

Let us first make the following observations, which can be proved
easily by induction: For any $k>0$ one has:\\
i) $(E,T^k(F_2))^i\ne0$ only for $i=k$ (see Remark \ref{rem:t2not}, i)).\\
ii) $(F_1,T^k(F_2))^k\ne0$ and $(F_1,T^k(F_2))^i=0$ for
$i>k$.\\
iii) $(F_1,T^k(F_2))^1\ne0$.\\
(Use the distinguished triangle $T^k(F_2)[k]\to E[2]\to
T^{k-1}(F_2)[k+1]$ and the non-trivial morphism $F_1\to E[2]$ for
ii) and $E[1-k]\to T^k(F_2)[1]\to T^{k+1}(F_2)[1]$ and
$(F_1,E[1-k])^1=0$ for iii).)

Hence, by ii) $(F_1,T^n(F_2))^n\ne0$, which by stability of $F_1$
and $T^n(F_2)$ shows
\begin{equation}\label{eqnii}
\phi(F_1)\leq\phi(T^n(F_2))+n\leq\phi(F_1)+2. \end{equation}
 On
the other hand, iii) implies
\begin{equation}\label{eqniii}
\phi(F_1)\leq\phi(T^n(F_2))+1\leq\phi(F_1)+2. \end{equation} This
leaves us with the possibilities $n=1,2,3$.

Next, $(E,F_1)^0\ne0$ and $(E,T^n(F_2))^n\ne0$ combined with Serre
duality yields
\begin{equation}\label{eqni}
\phi(E)<\phi(F_1)<\phi(E)+2~~{\rm
and}~~\phi(E)<\phi(T^n(F_2))+n<\phi(E)+2. \end{equation} (This
time equality can be excluded, as the spherical $E$ cannot be
isomorphic to a rigid $F_1[i]$ or $T^n(F_2)[i]$.)

If $n=3$, then $\phi(F_1)=\phi(T^n(F_2))+1$, which yields a
contradiction in (\ref{eqni}).

If $n=2$, then (\ref{eqnii}) and (\ref{eqniii}) show that
 $\phi(T^n(F_2))=\phi(F_1)$. (Use
that $F_1$, $F_2$, and $T^2(F_2)$ are numerically equivalent and
that therefore $\phi(T^2(F_2))-\phi(F_1)$ is an even integer.)
Both cases contradict (\ref{eqni}).

If $n=1$, we shall derive a contradiction as follows. For
notational simplicity assume $\phi(F_1)=0$, i.e.\
$Z(F_1)=Z(F_2)\in\RR_{>0}$. Now (\ref{eqnii}) reads
$-1\leq\phi(T(F_2))\leq 1$ and (\ref{eqni}) shows
$\phi(E)\in(-2,0)$. Since $E\to F_2\to T(F_2)$ is the
Harder--Narasimhan filtration of $F_2$, one  has
$0>\phi(E)\geq\phi(T(F_2))\geq -1$. The latter contradicts
$Z(E)+Z(T(F_2)) =Z(F_2)\in\RR_{>0}$.
\end{proof}

%%%%%%%%%%%%%%%%%%%%%%%%%%%%%%%%%%%%%%%%%%%%%%%%%%%%%%%%%%%%%%%%%%%%%%
\section{Generic twisted K3 and twisted abelian
surfaces}\label{sect:twist}

It is extremely difficult to obtain any information about the
space of stability conditions on a general triangulated category.
Even the most basic questions, e.g.\ whether it is non-empty or
connected, are usually very hard. In this section we aim at a
complete description of the space of stability conditions on  the
bounded derived category of coherent sheaves on a surface, which
is a generic twisted K3 surface or  an arbitrary twisted abelian
surface. Due to a technical problem, we are, for the time being,
only able to deal with the part of maximal dimension. In
particular, we will see that this space is connected and
simply-connected (see Theorem \ref{thm:generictwisted}).
Moreover, following ideas of Bridgeland, this can be used to
describe the group of autoequivalences in these cases (see
Theorem \ref{thm:genAut} and Corollary \ref{cor:twistder}).

\bigskip
%%%%%%%%%%%%%%%%%%%%%%%%%%%%%%%%%%%%%%%%%%%%%%%%%%%%%%%%%%%%%%%%%

\subsection{Stability conditions on twisted
surfaces}\label{subsect:twist} ~

\smallskip

For the bounded derived category of coherent sheaves on a
projective K3 surface Bridgeland describes in \cite{B} one
connected component $\Stab^\dag(X)$ of the space $\Stab(X)$ of
stability conditions on $\Db(X)=\Db(\coh(X))$ as a covering of a
certain period domain. In this section we will describe an
analogous connected component $\Stab^\dag(X,\alpha)$ of the space of
stability conditions on the bounded derived category
$\Db(X,\alpha)$ of $\alpha$-twisted coherent sheaves on a K3 or
abelian surface $X$ endowed with an additional Brauer class
$\alpha$. The construction follows closely the original one in
\cite{B} and we will therefore be brief and only explain the
necessary modifications.

\medskip

 Recall that the (cohomological) \emph{Brauer group} $\Br(X)$ of a
smooth complex projective variety $X$ is the torsion part of the
cohomology group $H^2(X,\ko_X^*)$ in the analytic topology (or,
equi\-valently, $H^2_{\mathrm{\acute{e}t}}(X,\ko_X^*)$ in the
tale topology). A \emph{twisted K3} (or \emph{abelian})
\emph{surface} is a pair $(X,\alpha)$ consisting of a K3
(abelian) surface $X$ and a Brauer class $\alpha\in\Br(X)$.

Any $\alpha\in\Br(X)$ can be represented by a \v{C}ech cocycle on
an open analytic cover $\{U_i\}_{i\in I}$ of $X$ using the
sections $\alpha_{ijk}\in\Gamma(U_i\cap U_j\cap
U_k,\mathcal{O}^*_X)$. An \emph{$\alpha$-twisted coherent sheaf}
$\kf$ consists of a collection $(\{\kf_i\}_{i\in
I},\{\varphi_{ij}\}_{i,j\in I})$, where $\kf_i$ is a coherent
sheaf on $U_i$ and $\varphi_{ij}:\kf_j|_{U_i\cap
U_j}\to\kf_i|_{U_i\cap U_j}$ is an isomorphism satisfying the
following conditions:
\[
\mbox{(1) $\varphi_{ii}=\mathrm{id}$;\;\;\;\;\;(2)
$\varphi_{ji}=\varphi_{ij}^{-1}$;\;\;\;\;\; (3)
$\varphi_{ij}\circ\varphi_{jk}\circ\varphi_{ki}=\alpha_{ijk}\cdot\mathrm{id}$.}
\]

By  $\coh(X,\alpha)$ we shall denote the abelian category of
$\alpha$-twisted coherent sheaves on $X$. Note that for different
\v{C}ech cocycles representing $\alpha$ the abelian categories
will be equivalent, although not canonical. For a discussion of
this see \cite{HS1,HS2}.

\begin{definition}
The bounded derived category  $\Db(\coh(X,\alpha))$ and the space
of numerical, locally finite stability conditions on it shall be
denoted
$$\Db(X,\alpha)~~\text{ resp. }~~\Stab(X,\alpha).$$
\end{definition}

In order to study the twisted categories $\Db(X,\alpha)$ and
equivalences between them, we have introduced in \cite{HS1} the
twisted Hodge structure $\widetilde H(X,\alpha,\ZZ)$ and the
twisted Chern character ${\rm ch}^\alpha:\Db(X,\alpha)\to
\widetilde H(X,\alpha,\ZZ)$. To `materialize' both
structures one needs to fix a B-field lift $B$ of the Brauer class
$\alpha$ and a cocycle representing $B$. So, the above twisted
Chern character ${\rm ch}^\alpha$ stands for ${\rm
ch}^B:\Db(X,\alpha_B)\to \widetilde H(X,B,\ZZ)$. We will largely
ignore this issue here, but the details, sometimes a little
confusing, can be found in \cite{HS2} (see also Remark
\ref{rem:expl}).

The twisted \emph{N\'eron--Severi group} is by definition
$${\rm NS}(X,\alpha):=\widetilde H^{1,1}(X,\alpha,\ZZ).$$
Note that for the trivial twist $\alpha=1$ the N\'eron--Severi
group ${\rm NS}(X,\alpha=1)$ differs from the classical
N\'eron--Severi group by  the additional sum $(H^0\oplus
H^4)(X,\ZZ)$. The Chern character ${\rm ch}^\alpha$ and the Mukai
vector $v^\alpha:={\rm ch}^\alpha\cdot (1,0,1)$ take values only
in ${\rm NS}(X,\alpha)$.

It is not difficult to show that ${\rm ch}^\alpha$ identifies the
numerical Grothendieck group with the N\'eron--Severi group, i.e.\
${\cal N}(\Db(X,\alpha))\cong{\rm NS}(X,\alpha)$.  In particular,
the stability function $Z$ of a numerical stability condition
 $\sigma=(Z,\kp)\in\Stab(\Db(X,\alpha))$ is of the form
$Z(E)=\langle v^\alpha(E),\varphi\rangle$
for some $\varphi\in{\rm NS}(X,\alpha)\otimes\CC$. This gives rise to the
\emph{period map}
$$\xymatrix{\pi:\Stab(X,\alpha)\ar[r]& {\rm NS}(X,\alpha)\otimes\CC,}~~
\xymatrix{\sigma=(Z,\kp)\ar@{|->}[r]&\varphi.}$$

\begin{definition} Denote
by $$P(X,\alpha)\subset{\rm NS}(X,\alpha)\otimes\CC$$  the open
subset of vectors $\varphi$ such that real part and imaginary part of
$\varphi$ generate a positive plane in ${\rm
NS}(X,\alpha)\otimes\RR$.
\end{definition}

\begin{remark}\label{rem:expl}
i) Suppose $B_0\in H^2(X,\QQ)$ is a B-field lift of $\alpha$,
i.e.\ $\alpha$ is the image of $B_0$ under the  exponential map
$H^2(X,\QQ)\to H^2(X,\ko^*_X)$.

To any real ample class $\omega\in H^{1,1}(X,\ZZ)\otimes\RR$ one
associates
$$\varphi:=\exp(B_0+i\omega)=1+(B_0+i\omega)+(B_0^2-\omega^2)/2+i(B_0.\omega)\in{\rm
NS}(X,\alpha)\otimes\CC.$$ Here we use that for a chosen B-field
lift $B_0$ of $\alpha$ the twisted cohomology $\widetilde
H^{1,1}(X,\alpha,\ZZ)$ can be identified with the integral part
of $\exp(B_0)\cdot\widetilde H^{1,1}(X,\QQ)$.

The real part $1+B_0+(B_0^2-\omega^2)/2$ and the imaginary part
$\omega+(\omega.B_0)$ of $\varphi$ span a plane which is positive
due to $\omega^2>0$.

As in the untwisted case, $P(X,\alpha)$ has two connected
components and we shall denote the one that contains
$\varphi=\exp(B_0+i\omega)$ as above by $P^+(X,\alpha)$. Thus, we
have
$$P^+(X,\alpha)\subset P(X,\alpha)\subset {\rm
NS}(X,\alpha)\otimes\CC$$ and one proves that the fundamental
group $\pi_1(P^+(X,\alpha))\cong\ZZ$ is generated by the loop
induced by the natural $\CC^*$-action
$(\lambda,\varphi)\mapsto\lambda\cdot\varphi$ on $P(X,\alpha)$.

ii) Note that if $B_0$ and $\omega$ are as in i), $B_1\in {\rm
NS}(X)\otimes\RR$, and $B:=B_1+B_0$, then
$$\varphi:=\exp(B+i\omega)=\exp(B_1)\cdot\exp(B_0+i\omega)\in{\rm
NS}(X,\alpha)\otimes\CC$$ is also contained in $P^+(X,\alpha)$.

iii) It is worth pointing out that in the twisted case there are
two sorts of B-fields. First, there are B-field lifts $B$ of the
Brauer class $\alpha$. It is the $(0,2)$-part of $B$ that matters
in this case. Second, as in the untwisted case, one needs B-fields
in order to `complexify' the polarization (or K\"ahler class) on
$X$. As a B-field lift of a given Brauer class $\alpha$ can always
be changed by a rational $(1,1)$-class without changing $\alpha$,
the difference between these two classes is not always clear cut.
\end{remark}

Following Bridgeland, we shall associate to any $\varphi\in
P^+(X,\alpha)$ a torsion theory and thus a $t$-structure on
$\Db(X,\alpha)$. Under additional conditions this will lead to a
stability condition, whose stability function is
\begin{equation}\label{eqn:slope}
Z_\varphi(E):=\left\langle v^\alpha(E),\varphi\right\rangle.
\end{equation}

Fix $\varphi\in P^+(X,\alpha)$ and let
$\kt,\kf\subset\coh(X,\alpha)$ be the following two full additive
subcategories: The non-trivial objects in $\kt$ are  the twisted
sheaves $E$ such that every non-trivial torsion free quotient
$E\twoheadrightarrow E'$ satisfies ${\rm Im}(Z_\varphi(E'))>0$. A
non-trivial twisted sheaf $E$ is an object in $\kf$ if $E$ is
torsion free and every non-zero subsheaf $E'$ satisfies ${\rm
Im}(Z_\varphi(E'))\leq0$.

It is easy to see that $(\kt,\kf)$ does define a torsion theory. In
particular, $\Hom(E,F)=0$ for $E\in\kt$ and $F\in\kf$ and any
sheaf $G\in\coh(X,\alpha)$ can be written in a unique way as an
extension
$$\xymatrix{0\ar[r]&E\ar[r]&G\ar[r]&F\ar[r]&0}$$ with
$E\in\kt$ and $F\in\kf$.

The heart of the induced $t$-structure is the abelian category
\[
\ka(\varphi):=\left\{E\in\Db(X,\alpha):\begin{array}{l}
\bullet\;\;\kh^i(E)=0\mbox{ for }i\not\in\{-1,0\},\\\bullet\;\;
\kh^{-1}(E)\in\kf,\\\bullet\;\;\kh^0(E)\in\kt\end{array}\right\}.
\]

Recall that a function $Z:\ka\to\CC$ on an abelian category is
called a \emph{stability function} if for all $0\ne E\in\ka$
either ${\rm Im}(Z(E))>0$ or $Z(E)\in\RR_{<0}$. The following is
the analogue of \cite[Lemma 6.2]{B}.

\begin{lem}\label{lem:slope}  For $\varphi=\exp(B+i\omega)$ as in Remark
\ref{rem:expl}, ii), the induced homomorphism
$$\xymatrix{Z_\varphi:\ka(\varphi)\ar[r]&\CC}$$ is a stability function on
$\ka(\varphi)$ if and only if for any spherical twisted sheaf
$E\in\coh(X,\alpha)$ one has $Z_\varphi(E)\not\in\RR_{\leq
0}$.\end{lem}

\begin{proof}
Let $E\in\ka(\varphi)$ with cohomology sheaves $\kh^i:=\kh^i(E)$,
$i=0,-1$. Then
$$Z_\varphi(E)=Z_\varphi(\kh^0)-Z_\varphi(\kh^{-1})=Z_\varphi(\kh^0_{\rm tor})
+Z_\varphi(\kh^0/\kh^0_{\rm tor})-Z_\varphi(\kh^{-1}),$$ where
$\kh^0_{\rm tor}$ denotes the torsion part of $\kh^0$.

The first thing to notice, is that for $\omega$ ample any torsion
sheaf $\ks\ne0$ satisfies ${\rm Im}(Z_\varphi(\ks))>0$ if
$\dim({\rm supp}(\ks))=1$ and $Z_\varphi(\ks)\in\RR_{<0}$ if
$\dim({\rm supp} (\ks))=0$. Secondly, by construction, ${\rm
Im}(Z_\varphi(\kh^0/\kh^0_{\rm tor}))>0$ if $\kh^0$ is not
torsion. Thus, $Z_\varphi$ is a stability function on
$\ka(\varphi)$ if and only if for any $E\in\kf$ with ${\rm
Im}(Z_\varphi(E))=0$ one has $Z_\varphi(E)\in\RR_{>0}$.

If one writes ${\rm ch}^\alpha(E)(1,0,1)=(r,\ell,s)$ with $r>0$,
then ${\rm Im}(Z_\varphi(E))=(\omega.\ell)-r(\omega.B)$ and $${\rm
Re}(Z_\varphi(E))=(1/2r)\left((\ell^2-2rs)+r^2\omega^2-(\ell-rB)^2\right).$$
Thus, ${\rm Im}(Z_\varphi(E))=0$ if and only if $(\ell-rB)$ is
orthogonal to $\omega$. The Hodge index theorem, applied to the
class $(\ell-rB)$, which really is of type $(1,1)$ on the
untwisted(!) surface, yields $(\ell-rB)^2\leq0$. On the other
hand, $\ell^2-2rs=-\chi(E,E)$. Using the existence of the usual
Jordan--H\"{o}lder filtration, one can restrict to the case that
$E$ is simple and hence $\ell^2-2rs\geq0$ or $=-2$ if $E$ is
spherical. This proves the assertion.\end{proof}

\begin{remark}
i) In fact, the proof shows that the condition is satisfied if
e.g.\ $\omega^2>2$.

ii) Later we shall be interested in the case when there are no
spherical objects. Then the lemma simply says that any
$\varphi=\exp(B+i\omega)$ defines a stability function. Here
$\omega$ varies in the K\"ahler cone which coincides with the
positive cone, since there would not be any $(-2)$-curves.
\end{remark}

The following two results are the twisted analogues of
Bridgeland's result in \cite{B}. For the first proposition we
assume that $\varphi=\exp(B+i\omega)$ is constructed as in Remark
\ref{rem:expl}, ii) and satisfies the condition of Lemma
\ref{lem:slope}, i.e.\ $Z_\varphi(E)\not\in\RR_{\leq0}$ for all
spherical sheaves.

\begin{prop}\label{prop:exist}  The stability function
 $Z_\varphi$ on the abelian category $\ka(\varphi)$
has the Harder--Narasimhan property and therefore defines a
stability condition $\sigma_\varphi$ on $\Db(X,\alpha)$ which,
moreover, is locally finite.\qqed
\end{prop}

Note that the proof of this result is rather round about already
in the untwisted case. For a direct proof when $B$ and $\omega$
are rational see \cite[Prop.\ 7.1]{B}. The argument works as well
in the twisted case.

\begin{definition} One denotes by
$\Stab^\dag(X,\alpha)$ the connected component of ${\rm Stab}(X,\alpha)$
that contains the stability conditions described by Proposition
\ref{prop:exist} (which form a connected set). Furthermore,
$U(X,\alpha)\subset\Stab^\dag(X,\alpha)$ shall denote the space of all
stability conditions in $\Stab^\dag(X,\alpha)$ for which all point
sheaves $k(x)$ are stable all of the same phase.

We say that a connected component of $\Stab(X,\alpha)$ is of
maximal dimension when the restriction of the period map
$\pi:\Stab(X,\alpha)\to\NS(X,\alpha)\otimes\CC$ to it is locally
homeomorphic.
\end{definition}

Clearly, $\Stab^\dag(X,\alpha)$ is a connected component of
maximal dimension. Also note that by definition $U(X,\alpha)$ is
connected.

\begin{prop}\label{prop:class}
Suppose $\sigma=(Z,\kp)\in\Stab(X,\alpha)$ is contained in a
connected component of maximal dimension, e.g.\ in
$\Stab^\dag(X,\alpha)$, and that for any closed point $x\in X$ the
skyscraper sheaf $k(x)$ is $\sigma$-stable of phase one with
$Z(k(x))=-1$. Then there exists $\varphi=\exp(B+i\omega)\in
P^+(X,\alpha)$ as in Remark \ref{rem:expl}, ii) such that the
heart of $\sigma$ coincides with $\ka(\varphi)$, i.e.\
$\kp((0,1])=\ka(\varphi)$.\qqed
\end{prop}

The result presumably continues to hold for any
$\sigma\in\Stab(X,\alpha)$. Note that the proposition does not
assert equality $\sigma=\sigma_\varphi$ of stability conditions,
but only of their hearts. The assumption that $\sigma$ is
contained in a connected component of maximal dimension is needed
in order to ensure the existence of generic deformations of
$\sigma$ and to eventually ensure that $\omega$ is really ample
(see Step 1 in the proof of \cite[Prop.\ 10.3]{B}).

\begin{remark}\label{rem:rationalbetter}
In fact, Bridgeland shows moreover that for `good' stability
conditions $\sigma$ one has $\sigma=\sigma_\varphi$ up to the
action of $\widetilde{\rm Gl}_2^+(\RR)$. For our purpose we need
the following stronger form of Proposition \ref{prop:class}:
Suppose $\sigma$ satisfies the conditions of the proposition and
its stability function $Z$ is of the form $Z_{\varphi'}$ with
$\varphi'\in{\rm NS}(X,\alpha)\otimes\QQ(i)$. Then $\varphi$ in
the assertion can be chosen such that up to the action of
$\widetilde{\rm Gl}_2^+(\RR)$ one has $\sigma_\varphi=\sigma$.

We sketch the argument in the untwisted case. Clearly, as
$Z(k(x))=-1$, one has $\varphi'=1+B'+i\omega'+(a+ib)$ with
$B',\omega'\in {\rm NS}(X)_\RR$ and $a,b\in H^4(X,\RR)$. Modulo
the action of $\widetilde{\rm Gl}_2^+(\RR)$ we can assume
$b=(B'.\omega')$. Indeed, if $\widetilde g$ acts on the stability
function $Z$ by
$g^{-1}=\left(\begin{array}{cc}1&\gamma\\0&1\end{array}\right)$
with $\gamma=((B'.\omega')-b)/(\omega.\omega)$, then $g^{-1}Z$
has this property. Once this is achieved, Bridgeland shows that
$\omega'$ is ample, so he sets $\omega=\omega'$, and that one can
choose $B=B'$. The difference between the two stability functions
$Z_{\varphi'}$ and $Z_\varphi$, where $\varphi=\exp(B+i\omega)$,
is the real degree four part.

Now one uses the additional assumption $B'$, $\omega'$, $a$, and
$b$ are all rational. The special $g$ used to modify $\varphi'$
does not affect this. One first checks that then $a<(B.B)/2$. In
order to see this, pick $r>0$ such that $rB$ is integral and such
that $r^2(B.B)$ is divisible by $2r$. Then there exists a
$\mu$-stable vector bundle $E$ with ${\rm c}_1(E)=rB$, such that
$\langle v(E),v(E)\rangle=0$. Then $E[1]\in \ka(\varphi)$ and
${\rm Im}(Z_{\varphi'}(E))=0$. On the other hand, ${\rm
Re}(Z_{\varphi'}(E))>0$ if and only if $a<(B.B)/2$.

Eventually apply an element  of $\widetilde{\rm Gl}_2^+(\RR)$
which acts by
$g^{-1}=\left(\begin{array}{cc}1&0\\0&(B.B)-2a\end{array}\right)$
on $Z_{\varphi'}$. The resulting stability conditions is of the
form $Z_\varphi$ with $\varphi=\exp(B+i((B.B)-2a)\omega)$.
\end{remark}

\begin{prop}
The period map defines a covering map $\pi:\Stab^\dag(X,\alpha)\to
P^+_0(X,\alpha)$ onto some open subset $P^+_0(X,\alpha)\subset
P^+(X,\alpha)$.\qqed
\end{prop}

Moreover, the subgroup of autoequivalences that preserve the
component $\Stab^\dag(X,\alpha)$ acts freely on
$\Stab^\dag(X,\alpha)$.

Of course, as in the untwisted case $P^+_0(X,\alpha)$ is cut out
by all hyperplanes orthogonal to $(-2)$-classes, but we will not
need this.

Let $\sigma\in U(X,\alpha)$. Then $\sigma$ can be changed by an
element in $\widetilde{\rm Gl}_2^+(\RR)$ such that all point
sheaves $k(x)$ are stable of phase one with $Z(k(x))=-1$.

\begin{remark}
In complete analogy to the discussion in Sections 11, 12, and 13
in \cite{B} one also proves the following result: If $(X,\alpha)$
is a twisted K3 or abelian surface such that $\coh(X,\alpha)$ does
not contain any spherical objects, then
$U(X,\alpha)=\Stab^\dag(X,\alpha)$.

However, a stronger result with a more direct proof relying only
on the discussion of Section \ref{sect:GenCat} will be established
further below (see Theorem \ref{thm:generictwisted}).

For untwisted abelian surfaces, which never admit any spherical
sheaves, this was remarked in \cite[Sect.\ 15]{B}.
\end{remark}

\bigskip

%%%%%%%%%%%%%%%%%%%%%%%%%%%%%%%%%%%%%%%%%%%%%%%%%%%%%%%%%%%%%%%%%%%
\subsection{Twisted surfaces without spherical
sheaves}\label{subsect:twistnosph} ~

\smallskip

Typical examples of spherical sheaves on a K3 surface $X$ are the
trivial line bundle $\ko_X$ and the structure sheaves $\ko_C$ of
$(-2)$-curves $\PP^1\cong C\subset X$. The latter do not exist on
a generic projective K3 surface, but e.g.\ $\ko_X$, of course,
persists. The situation is different in the twisted case.

\begin{definition}\label{def:generic}
A twisted K3 or abelian surface $(X,\alpha)$ will be called
\emph{generic} if $\coh(X,\alpha)$ (or, equivalently,
$\Db(X,\alpha)$) does not contain any spherical objects.
\end{definition}

\begin{remark}\label{rem:collect}
Here we collect a few observations, that will be useful
throughout this section and in particular explain why in the
above definition the two formulations are really equivalent. All
the assertions can be proved by adapting the arguments of Section
\ref{subsect:sphercial}. One replaces stability conditions on
triangulated categories by ordinary (Gieseker) stability of
twisted coherent sheaves (see \cite{Y}). So, in the following,
semistable and stable factors of a sheaf are meant with respect
to  the usual Harder--Narasimhan respectively Jordan--H\"older
filtrations in the context of Gieseker stability.

 i) Let us first work in the abelian category $\coh(X,\alpha)$ (cf.\ Proposition
 \ref{prop:stablefact}):\\
 $\bullet$ The  stable factors of a rigid sheaf $E\in\coh(X,\alpha)$
 are spherical.\\
 $\bullet$ At most one stable factor of a semi-rigid sheaf $E\in\coh(X,\alpha)$
 is semi-rigid and all others are spherical.

 In particular, if there are no spherical sheaves in
 $\coh(X,\alpha)$, then there are no rigid sheaves and semi-rigid
 ones are stable (cf.\  Corollary \ref{lem:skyscap}).

 ii) In the derived category one has the following:\\
 $\bullet$ If $E\in\Db(X,\alpha)$ is rigid, then all cohomology
 sheaves $\kh^i(E)$ are rigid. Thus, if there are no spherical sheaves
 in $\coh(X,\alpha)$, then $\Db(X,\alpha)$ does not contain any
 rigid or spherical objects. To prove this one applies
 \cite[Lemma 2.9]{BM}, valid also in the twisted context, to show
$\sum(\kh^i(E),\kh^i(E))^1\leq (E,E)^1=0.$
 \\
 $\bullet$ Similarly, if $E\in\Db(X,\alpha)$ is semi-rigid, then at most one
 cohomology sheaf $\kh^i(E)$ is semi-rigid and all others are
 rigid. Thus, if there are no spherical objects in
 $\coh(X,\alpha)$, a semi-rigid object in $\Db(X,\alpha)$ is up to
 shift isomorphic to a semi-rigid, stable sheaf.
\end{remark}

Twisted abelian surfaces are always generic in this sense, as they
never admit any spherical sheaves. (Note that the twist $\alpha$
in a twisted surface $(X,\alpha)$ might very well be trivial,
i.e.\ $(X,\alpha)$ might simply be the surface $X$. Thus,
ordinary abelian surfaces are covered by our discussion.)

As we shall see in Section \ref{sect:densi} generic twisted K3
surfaces are dense  among all twisted K3 surfaces. Thus, the
following consequence of Proposition \ref{prop:stablefact} (see
also Corollary \ref{lem:skyscap}) applies to a dense subset of
twisted K3 surfaces and all twisted abelian surfaces.

\begin{cor}\label{cor:ptsarestable}
Suppose $(X,\alpha)$ is a generic twisted K3 surface or an
arbitrary twisted abelian surface. Then the skyscraper sheaves
$k(x)$ for closed points $x\in X$ are $\sigma$-stable for all
$\sigma\in\Stab(X,\alpha)$.
\end{cor}

\begin{proof}
Indeed, for a closed point the skyscraper sheaf $k(x)$ is
semi-rigid.
\end{proof}

This then yields

\begin{thm}\label{thm:generictwisted}
Let $(X,\alpha)$ be a generic twisted K3 surface or an arbitrary
twisted abelian surface. Then $\Stab(X,\alpha)$ admits only one
connected component of maximal dimension.

Moreover, this maximal component is simply-connected and the
restriction of the  period map
$\xymatrix{\pi:\Stab(X,\alpha)\ar[r]& {\rm
NS}(X,\alpha)\otimes\CC}$ to it can be viewed as the universal
cover of $P^+(X,\alpha)$.
\end{thm}

\begin{proof}
Let $\sigma\in\Stab(X,\alpha)$ be a stability condition contained
in a connected component of maximal dimension. In particular,
after a small deformation the period of $\sigma$ will be rational,
i.e.\ contained in ${\rm NS}(X,\alpha)\otimes\QQ(i)$.

Due to Corollary \ref{cor:ptsarestable} all point sheaves $k(x)$
are stable.  Now choose a  torsion free semi-rigid twisted sheaf
$E$ on $(X,\alpha)$. The existence of such an $E$ can be deduced
from standard existence results, see e.g.\ \cite{Y} for the
twisted case. Note that Corollary \ref{lem:skyscap} also applies
to $E$ and shows that $E$ is stable with respect to $\sigma$.

Since $E$ is torsion free, there exists a non-trivial $E\to k(x)$
for any closed point $x\in X$. Hence, the assumptions of Lemma
\ref{prop:phase} are satisfied for any two points $x_1\ne x_2\in
X$. Thus, either $\phi(k(x_1))=\phi(k(x_2))$ or $k(x_1)\cong
k(x_2)[\pm2]$, but the latter is absurd. (For the equality of the
phases of the point sheaves $k(x)$ see also Remark
\ref{rem:avoid}, which explains how to avoid to work with $E$.)
Together with Remark \ref{rem:rationalbetter}, this shows that
$\sigma\in U(X,\alpha)$. In particular,
$U(X,\alpha)=\Stab^\dag(X,\alpha)$ and $\Stab^\dag(X,\alpha)$ is
the only connected component of maximal dimension.

Next, consider the period map $\pi:U(X,\alpha)\to P^+(X,\alpha)$.
As mentioned before, the fundamental group of $P^+(X,\alpha)$ is
the free cyclic group generated by the loop produced by the
$\CC^*$-action which lifts to the double shift
$\sigma\mapsto\sigma[2]$ on $\Stab(X,\alpha)$. Thus, in order to
prove the remaining assertions, it suffices to show that
$\pi:U(X,\alpha)\to P^+(X,\alpha)$ is surjective.

 As there are no spherical objects in $\coh(X,\alpha)$, the surface $X$ does not,
 in particular, contain any $(-2)$-curve. Thus, K\"ahler cone and positive cone coincide.
 Therefore, any element in $P^+(X,\alpha)$ is up to the action
${\rm Gl}_2^+(\RR)$ of the form
 $\varphi$ with $\varphi= \exp(B+i\omega)$
 as in Remark \ref{rem:expl}, ii). Hence, due to Lemma \ref{lem:slope} and
 Proposition \ref{prop:exist} all the periods are in the image of
 $\pi:U(X,\alpha)\to P^+(X,\alpha)$.
 \end{proof}

\begin{remark}\label{rem:excuse}
The result can be strengthened as follows. All stability
conditions which can be deformed to a stability condition with
rational period form one connected component. This component,
moreover, is simply-connected and the universal cover of
$P^+(X,\alpha)$.

Moreover, following \cite{B}, one shows that if $\sigma\in\Stab(X,\alpha)$ and
$\pi(\sigma)\in\kp^+(X,\alpha)$ then $\sigma$ is in the unique
component $\Stab^\dag(X,\alpha)$ of maximal dimension.
\end{remark}

\bigskip

%%%%%%%%%%%%%%%%%%%%%%%%%%%%%%%%%%%%%%%%%%%%%%%%%%%%%%%%%%%%
\subsection{(Auto)equivalences}\label{subsect:twisauto}
~

\smallskip

As autoequivalences of the derived category $\Db(X,\alpha)$ act
naturally on the space of stability conditions $\Stab(X,\alpha)$,
the results of the previous section can be used to fully determine
the group of autoequivalences $\Aut(\Db(X,\alpha))$ for a generic
twisted K3 surface $(X,\alpha)$. Everything is known for
arbitrary (twisted) abelian surfaces, so we will concentrate on
K3 surfaces here. In particular, we shall show that
autoequivalences of $\Db(X,\alpha)$ induce orientation preserving
Hodge isometries of $\widetilde H(X,\alpha,\ZZ)$, i.e.\ that the
orientation of the four-space of  positive directions is not
changed. This had first been conjectured by Szendr\H{o}i, but
even for untwisted K3 surfaces it is still not known in general
(see e.g.\ \cite{HFM}).

In the following we shall tacitly use that any equivalence is of
Fourier--Mukai type, which was proved by Orlov in \cite{Or1} in
the untwisted case and in \cite{CS} in general.

\begin{thm}\label{thm:genAut}
Suppose $(X,\alpha)$ is a generic twisted K3 surface. Then there
exists a  short exact sequence
$$\xymatrix{0\ar[r]&\ZZ\ar[r]&\Aut(\Db(X,\alpha))\ar[r]&\Aut^+(\widetilde H(X,\alpha,\ZZ))\ar[r]&0,}$$
where $\Aut^+(\widetilde H(X,\alpha,\ZZ))$ denotes the group of
orientation preserving Hodge isometries.
\end{thm}

\begin{proof}
Let ${\rm id}\ne\Phi\in\Aut(\Db(X,\alpha))$. Then the
cohomological Fourier--Mukai transform defined in terms of the
twisted Mukai vector of the Fourier--Mukai kernel of $\Phi$
yields a Hodge isometry $\Phi_*\in\Aut(\widetilde
H(X,\alpha,\ZZ))$ (see \cite{HS1}). Clearly, $\Phi_*$ preserves
the orientation of the positive plane $(\widetilde
H^{2,0}\oplus\widetilde H^{0,2})(X,\RR)$. Thus, in order to show
that it preserves the orientation of all four positive
directions, it suffices to prove that the orientation of the two
remaining positive directions in ${\rm NS}(X,\alpha)$ is
preserved.

If $\sigma\in\Stab(X,\alpha)$ is the stability condition
associated to $\varphi=\exp(B+i\omega)$ as in Proposition
\ref{prop:exist}, then the stability function of $\Phi(\sigma)$ is
$Z_{\Phi_*\varphi}$. Since there is only one connected component of
maximal dimension in $\Stab(X,\alpha)$  and $\pi$ is continuous,
$\pi(\sigma)$ and $\pi(\Phi(\sigma))$ are contained in the same
connected component of $P(X,\alpha)$, i.e.\ in $P^+(X,\alpha)$.
In other words, real and imaginary part of $\varphi$ and of
$\Phi_*(\varphi)$ yield the same orientation of the positive
directions. Hence $\Phi_*\in\Aut^+(\widetilde H(X,\alpha,\ZZ))$.

If $\Phi_*={\rm id}$, then for any stability condition $\sigma$
the induced stability condition $\Phi(\sigma)$ has the same
stability function. Thus, $\Phi$ acts as a non-trivial deck
transformation of the covering map $\pi:U(X,\alpha)\to
P^+(X,\alpha)$. As a generator of $\pi_1(P^+(X,\alpha))\cong\ZZ$
lifts to the shift $E\mapsto E[2]$ acting as a deck
transformation of the simply-connected space $U(X,\alpha)$, one
finds that $\Phi$ is an even shift $E\mapsto E[2n]$.
\end{proof}

The fact that connectivity  of the space of stability conditions
can be used to prove that autoequivalences preserve the
orientation of the positive directions in $\widetilde H(X)$ was
observed by Bridgeland \cite[Conj.\ 1.2]{B}. In the following we
shall give another proof of this fact which is more direct and
does not rely on the concept of a stability condition and, in
particular, not on Theorem \ref{thm:generictwisted}.

\begin{prop}\label{prop:shiftsheaf} Let $(X,\alpha)$ and $(Y,\beta)$ be
generic twisted K3 surfaces. Suppose
$$\xymatrix{\Phi_\ke:\Db(X,\alpha)\ar[r]^-\sim&\Db(Y,\beta)}$$ is a Fourier--Mukai
equivalence with $\ke\in\Db(X\times
Y,\alpha^{-1}\boxtimes\beta)$. Then there exists a twisted
sheaf(!) $\kf\in\coh(X\times Y,\alpha^{-1}\boxtimes\beta)$ such
that $\ke\cong\kf[n]$, for some $n\in\ZZ$.\end{prop}

\begin{proof} For $x\in X$ a closed point, let
$\ke_x:=\ke|_{\{x\}\times Y}=\Phi_\ke(k(x))$. Since $\Phi_\ke$ is
an equivalence, $\ke_x$ is a semi-rigid object in $\Db(Y,\beta)$.

By Remark \ref{rem:collect}, ii)  there exists at most one
$i\in\ZZ$ such that $\kh^i(\ke_x)\neq 0$.  It is not difficult to
see that $i$ does not depend on the point $x\in X$. This proves
that there exists a sheaf $\kf\in\coh(X\times
Y,\alpha^{-1}\boxtimes\beta)$ and an integer $n$ such that
$\ke=\kf[n]$.
\end{proof}

\begin{cor}\label{cor:orient} Suppose $(X,\alpha)$ and $(Y,\beta)$ are
generic twisted K3 surfaces.

Then any Fourier--Mukai equivalence
$\Phi_\ke:\Db(X,\alpha)\congpf\Db(Y,\beta)$ induces an
orientation preserving Hodge isometry $\Phi_{\ke*}:\widetilde
H(X,\alpha,\ZZ)\congpf\widetilde H(Y,\beta,\ZZ)$.\end{cor}

\begin{proof} Note  again that in order to define the cohomological
twisted Fourier--Mukai functor $\Phi_{\ke*}$ one needs to fix
B-field lifts of $\alpha$ and $\beta$, but we will suppress all
technical details here.

By Proposition \ref{prop:shiftsheaf} the kernel $\ke$ is, up to
shift, a flat family of semi-rigid and hence, by Remark
\ref{rem:collect}, ii), stable sheaves on $Y$ parametrized by
$X$. Now one argues as in \cite[Prop.\ 5.5, Rem.\ 5.7]{HS1} in
order to deduce that the Hodge isometry $\Phi_{\ke*}:\widetilde
H(X,\alpha,\ZZ)\congpf\widetilde H(Y,\beta,\ZZ)$ is orientation
preserving.
\end{proof}

\begin{remark} Due to the previous proposition, there are twisted K3
surfaces $(X,\alpha)$ such that the Hodge isometry
$$j:=\mathrm{id}_{H^0\oplus H^4}\oplus-\mathrm{id}_{H^2}:\widetilde
H(X,\alpha,\ZZ)\rightarrow\widetilde H(X,\alpha^{-1},\ZZ)$$ is not
induced by any (Fourier--Mukai) equivalence.\end{remark}

In \cite[Thm.\ 0.1]{HS2} we proved that any orientation preserving
Hodge iso\-metry $\widetilde H(X,\alpha,\ZZ)\cong \widetilde
H(Y,\beta,\ZZ)$ between arbitrary  twisted K3 surfaces can be
lifted to a derived equivalence.

Thus, as a consequence of Corollary \ref{cor:orient} and of
\cite[Thm.\ 0.1]{HS2},  we obtain the following more precise
result, which is expected to hold without any genericity
assumption.

\begin{cor}\label{cor:twistder} Let $(X,\alpha)$ and $(Y,\beta)$ be
two generic twisted K3 surfaces. Then a Hodge isometry
$$g:\widetilde H(X,\alpha,\ZZ)\congpf\widetilde H(Y,\beta,\ZZ)$$
 can be lifted to a Fourier--Mukai equivalence
$$\Db(X,\alpha)\congpf\Db(Y,\beta)$$ if and only if $g$ is
orien\-tation preserving.\qqed
\end{cor}

\bigskip

%%%%%%%%%%%%%%%%%%%%%%%%%%%%%%%%%%%%%%%%%%%%%%%%%%%%%%%%%%%%%%%%%%%%
\subsection{Density of generic twisted K3
surfaces}\label{sect:densi} ~

\smallskip

Here we shall explain that generic twisted K3 surfaces in the
sense of Definition \ref{def:generic} are dense in the moduli
space. Unfortunately, for technical reasons this information seems
difficult to use in order to pass by deformation from the generic
case to the case of an arbitrary K3 surface.

\begin{lem}\label{lem:nospher} For any K3 surface $X$ with Picard
number one, there exist infinitely many $\alpha\in\Br(X)$ such
that $(X,\alpha)$ is generic, i.e.\ $\coh(X,\alpha)$ does not
contain any spherical objects.\end{lem}

\begin{proof} We shall work with specific B-field lifts $B$
of the Brauer class $\alpha=\alpha_B$ and with the Mukai vector
$v^B={\rm ch}^B\cdot(1,0,1)$. For the notation compare
\cite{HS1,HS2}.

Let $\Pic(X)\cong\ZZ H$ and suppose $H^2=2d>2$. Consider B-fields
$B\in H^2(X,\QQ)$ such that: \begin{itemize}\item[(1)] $(B.H)=(B.
B)=0$ and\item[(2)] There exists a prime number $p$ dividing $d$
and the order of $\alpha_B$.
\end{itemize}

Under these hypotheses, the quadratic form associated to
$\NS(X,\alpha_B)$ is represented by the matrix
\[
M=\left(\begin{array}{ccc}
H^2 & 0 & 0\\
0 & 0 & mp\\
0 & mp & 0
\end{array}\right),
\]
for some integer $m$. Of course, $M$ does not represent $-2$,
which excludes the existence of spherical objects
$E\in\Db(X,\alpha_B)$, for a spherical $E$ satisfies $\langle
v^B(E),v^B(E)\rangle=-2$.

If $H^2=2$, then consider B-fields $B$ satisfying (1) and such
that $-1$ is not a square modulo the order of $\alpha_B$. The
intersection form on $\NS(X,\alpha_B)$ has matrix form $M$ and
does not represent $-2$.\end{proof}

Let $\Lambda:=U^{\oplus 3}\oplus (-E_8)^{\oplus 2}$ and
$\widetilde\Lambda:=U^{\oplus 4}\oplus (-E_8)^{\oplus 2}$. The
space $Q:=\{x\in\PP(\Lambda\otimes\CC):x^2=0,(x,\bar x)>0\}$ is
the period domain of all K3 surfaces while $Q_{\rm alg}\subset Q$
is the dense subset of periods of algebraic K3 surfaces. In other
words, $Q_{\rm alg}$ is the set of those $x\in Q$ for which there
exists a class $h\in x^\perp\cap\Lambda$ with $h^2>0$. According
to \cite{Hu}, the space $\widetilde
Q:=\exp(\Lambda\otimes\RR)(Q_\mathrm{alg})\subset
\PP(\widetilde\Lambda\otimes\CC)$ is the period domain of the
moduli space of generalized Calabi--Yau structures on algebraic K3
surfaces obtained by (real) B-field shifts. Using Lemma
\ref{lem:nospher} it is very easy to derive the following result.

\begin{cor}\label{cor:nospher} The periods of generic twisted K3 surfaces
$(X,\alpha)$ are dense in $\widetilde Q$.\qqed\end{cor}

More precisely, the corollary says that the periods of K3 surfaces
$X$ shifted by B-fields $B$ inducing a generic twisted K3 surface
$(X,\alpha_B)$ form a  dense subset of $\widetilde Q$.
%%%%%%%%%%%%%%%%%%%%%%%%%%%%%%%%%%%%%%%%%%%%%%%%%%%%%%%%%%%%%%%%%%%
\section{Generic non-algebraic K3 surfaces}\label{sect:nonalgK3}

In this section we shall deal with generic non-algebraic K3
surfaces. The abelian category $\coh(X)$ and hence its derived
category $\Db(X)$ of a non-algebraic K3 surface is smaller than
in the projective case and contrary to Gabriel's classical result
for algebraic varieties $\coh(X)$ usually does not determine $X$
(see \cite{V4}).

However, although the generic non-algebraic K3 surface $X$ has
trivial Picard group, there always is the trivial line bundle
$\ko_X$ which gives rise to a spherical object in $\coh(X)$. In
this sense, the generic non-algebraic case is `less generic' than
the generic twisted projective case discussed in the previous
section. This makes a study of the generic non-algebraic situation
interesting; it leads to K3 categories with exactly one spherical
object (up to shift) and the method of Section \ref{sect:onlyone}
will be applicable.

\bigskip
%%%%%%%%%%%%%%%%%%%%%%%%%%%%%%%%%%%%%%%%%%%%%%%%%%%%%%%%%%%%%%%%%%%%%

\subsection{Spherical objects on generic K3 surfaces}\label{sect:gennonpro} ~

\smallskip

We will henceforth call a K3 surface $X$ \emph{generic} if its
Picard group is trivial. Clearly, with this definition a generic K3
surface is never algebraic and $H^{1,1}(X,\ZZ)=0$. The latter is
equivalent to $\kn(X)\cong (H^0\oplus H^4)(X,\ZZ)$.

Notice that, in principle, when dealing with smooth compact analytic
(non-projective) va\-rie\-ties $X$, we should work with the full
subcategory $\Db_{\rm coh}(\Mod{X})$ of $\Db(\Mod{X})$ consisting of
complexes of $\ko_X$-modules with coherent cohomology. Indeed due to
the general results in \cite{Spalt} the usual derived functors
(e.g.\ push-forwards, tensor products and pull-backs) are
well-defined in this category and it makes perfect sense to
introduce the notion of Fourier-Mukai functors. Moreover Serre
duality holds true in $\Db_{\rm coh}(\Mod{X})$ (see, for example,
\cite[Prop.\ 5.1.1]{BB}). If $X$ is now a smooth compact analytic
surface, then the natural functor $\Db(X)\to\Db(\Mod{X})$ induces an
equivalence $\Db(X)\iso\Db_{\rm coh}(\Mod{X})$ (see for example
\cite[Prop.\ 5.2.1]{BB}) and any $E\in\coh(X)$ admits a locally free
resolution of finite length (see \cite{Sc}). Since for the rest of
this paper we will just consider generic K3 surfaces, we will continue to work
with the triangulated category $\Db(X)$.

\begin{lem}
If $X$ is a generic K3 surface, then the trivial line bundle
$\ko_X$ is up to shift the only spherical object in $\Db(X)$.
\end{lem}

\begin{proof}
Suppose $E\in\coh(X)$ is a rigid sheaf and let $v(E)=(r,0,s)$ be
its Mukai vector. Then $2\leq\chi(E,E)=2rs$ and since $r\geq0$,
one finds $r,s>0$. In particular, $E$ cannot be torsion.

Next, we prove that a rigid $E$ is torsion free. To this end,
consider the short exact sequence  $0\to E_{\rm tor}\to E\to
E'\to 0$ with $E'$ torsion free and $E_{\rm tor}$ the torsion
subsheaf of $E$. In particular, $\Hom(E_{\rm tor},E')=0$ and
\cite[Cor.\ 2.8]{Mu} applies (the abelian analogue of Lemma
\ref{lem:Mukaigen}). Thus, $E_{\rm tor}$ is rigid as well, which
contradicts the above observation.

Again following Mukai, one proves that a rigid $E$ is in fact
locally free. Indeed, if not, then the surjection $E^{\vee\vee}\to
E^{\vee\vee}/E$ from the locally free reflexive hull
$E^{\vee\vee}$ to the torsion sheaf $E^{\vee\vee}/E\ne0$
concentrated in dimension zero can be deformed providing
non-trivial deformations of $E$ which is excluded by rigidity.

We shall prove that any rigid sheaf is isomorphic to
$\ko_X^{\oplus r}$ for some $r$. This will be done by induction
on the rank. The case $r=1$ is obvious.

Observe that $\chi(\ko_X,E)=r+s>0$. Hence, $\Hom(\ko_X,E)\ne0$ or
$\Hom(E,\ko_X)\ne0$, but it suffices to consider the first  case.
Indeed if $\Hom(E,\ko_X)\ne0$, then $\Hom(\ko_X,E\dual)\ne0$ and
$E\dual$ is still rigid.  Now suppose $\ko_X\to E$ is non-trivial
and hence  injective, for $E$ is torsion free. By induction we
may assume that we have a short exact sequence $0\to\ko_X^{\oplus
r}\to E\to F\to0$ with $r>0$. We shall show that
$\Hom(\ko_X,F)=0$ if $r$ is chosen maximal.

We first prove that $F$ is torsion free as well. Consider the
natural short exact sequence $0\to \ko_X^{\oplus r}\to G\to
F_{\rm tor}\to0$, where $G$ is the kernel of the surjection
$E\twoheadrightarrow F \twoheadrightarrow F/F_{\rm tor}$. Since
$\Ext^1(F_{\rm tor},\ko_X)=0$ and $E$ is torsion free, this
yields $F_{\rm tor}=0$. Suppose now that $\Hom(\ko_X,F)\ne0$.
Using $\Ext^1(\ko_X,\ko_X)=0$, any $\ko_X\to F$ can be lifted to
a morphism $\ko_X\to E$. Since $F$ is torsion free, one thus
obtains an injection $\ko_X^{\oplus r+1}\hookrightarrow E$.
Clearly, this process will terminate.

If $\Hom(\ko_X,F)=0$, then we are once more in the situation of
\cite[Cor.\ 2.8]{Mu}. Hence $F$ is rigid of smaller rank and thus
$F\cong\ko_X^{\oplus r'}$ and therefore $E\cong\ko_X^{\oplus
r+r'}$.

Proposition \ref{rem:mod} then concludes the proof.
\end{proof}

\bigskip
%%%%%%%%%%%%%%%%%%%%%%%%%%%%%%%%%%%%%%%%%%%%%%%%%%%%%%%%%%%%%%%%%%%%%

\subsection{Construction of stability conditions on
generic K3 surfaces} ~

\smallskip

Instead of adapting Bridgeland's discussion to the non-projective
case, we shall give an ad hoc approach that works well for generic
K3 surfaces. So, in the following, $X$ will denote a K3 surface
with trivial Picard group. Throughout we will use that torsion
sheaves on such a K3 surface are concentrated in points and that
the reflexive hull of an arbitrary torsion free sheaf is a
$\mu$-semistable vector bundle with trivial determinant.

The first step consists of an explicit construction of
certain stability conditions. This immitates arguments
in \cite{B} but differs at a few places. In particular,
the actual check of all the conditions is more direct and
in some cases also the heart has a different shape.

Consider the open subset
\begin{equation}\label{decompR}R:=\CC\setminus \RR_{\geq-1}=R_+\cup R_-\cup R_0,\end{equation}
where
$$R_+:=R\cap\HH,~~R_-:=R\cap(-\HH),~~\text{and}~~R_0:=R\cap\RR$$
with $\HH$ denoting the upper half-plane.
%\begin{itemize}
%\item[{\bf i)}] $R_+$ is the set of complex numbers $z=x+iy\in R$ with $y>0$;
%\item[{\bf ii)}] $R_0:=(-\infty,-1)$;
%\item[{\bf iii)}]  $R_-$ is the set of complex numbers $z=x+iy\in R$ with $y<0$.
%\end{itemize}

 For any $z=u+iv\in R$  we consider a torsion theory
$$\kf(z),\kt(z)\subset\coh(X)$$ which is defined uniformly as
follows: We let $\kf(z)$ be the full subcategory of all torsion
free sheaves of degree $\leq v$, whereas $\kt(z)$ is the full
subcategory that contains all torsion sheaves and all torsion
free sheaves of degree $>v$. The degree is taken with respect to
any K\"ahler structure, but this does not matter at all. In fact,
if $z\in R_+\cup R_0$ then $\kf(z)$ and $\kt(z)$ are simply the
full subcategories of all torsion free sheaves respectively
torsion sheaves and if $z\in R_-$ then $\kf(z)$ is trivial and
$\kt(z)=\coh(X)$.

Consequently, the tilt $\ka(z)$ of $\coh(X)$ with respect to the torsion
theory $(\kf(z),\kt(z))$ yields the abelian category (the heart of the
corresponding $t$-structure):
$$\ka:=\ka(z)=\left\{E~|~\kh^{-1}(E)~\mbox{torsion free},~\kh^0(E)~\mbox{torsion},\kh^i(E)=0
 \mbox{ for }i\ne0,-1\right\}$$ for $z\in R_+\cup R_0$ and
$$\ka(z)=\coh(X)$$ for $z\in R_-$.

For every $z=u+iv\in R$ one defines the function
$$\xymatrix{Z:\ka(z)\ar[r]&\CC},~
\xymatrix{E\ar@{|->}[r]&\langle v(E),(1,0,z)\rangle=-u\cdot
r-s-i(r\cdot v).}$$ Here, as usual $v(E)=(r,0,s=r-{\rm c}_2(E))$.

\begin{lem}\label{lem:stabfctgen}
For any $z\in R$ the function $Z$ defines a stability function
on $\ka(z)$.
\end{lem}

\begin{proof} We have to show that for $0\ne E\in\ka(z)$ one
has $Z(E)\in\HH\cup \RR_{<0}$.

For $z\in R_+$ one has ${\rm Im}(Z(E[1]))>0$ for any non-trivial
torsion free sheaf $E$. If $E$ is a torsion sheaf, i.e.\
$E$ is a finite length sheaf concentrated in dimension zero,
then $Z(E)={\rm c}_2(E)=-\chi(E)<0$.

For $z\in R_0$ torsion sheaves can be dealt with in the same way.
For a torsion free sheaf $E$, one has this time ${\rm
Im}(Z(E))=0$. But the real part of $Z(E)$ is $-u\cdot
r-s=-(u+1)\cdot r+{\rm c}_2(E)$. Using the additivity of the Mukai vector for short exact sequences one can easily reduce to the case that neither $E$ nor its dual contains $\ko_X$. At that point one applies $0\geq\chi(E)=-{\rm c}_2(E)+2\rk(E)$, which yields
the even stronger inequality ${\rm c}_2(E)\geq2\rk(E)$.

If $z\in R_-$, then ${\rm Im}(Z(E))=-r\cdot v>0$ for any
non-torsion sheaf $E$ and for torsion sheaves one concludes as
before.
\end{proof}

\begin{remark}\label{rem:onlysemirig}
Let us next pass to the classification of minimal and stable
objects in $\ka(z)$.

i) For $z\in R_0$ every object in $\ka(z)=\ka$ is semistable. For
$z\in R_-$ an object $E\in \ka(z)=\coh(X)$ is semistable if $E$ is
torsion  or torsion free. Finally, for $z\in R_+$ torsion sheaves
and shifted vector bundles $E[1]$ define semistable objects in
$\ka(z)=\ka$. There are, however, other semistable objects. As an
example one can consider for any closed point $x\in X$ the unique
non-trivial extension $0\to\ko_X[1]\to F\to k(x)\to 0$.

ii) It is well-known that the minimal objects of $\coh(X)$ are
the point sheaves $k(x)$. This covers the case $z\in R_-$.

For $z=u+iv\in R_+\cup R_0$, i.e.\ for $v\geq0$, the minimal
objects of the category are classified as follows (see e.g.\ the
discussion in \cite{H}): They consist of the point sheaves
$k(x)$ and the shifted $\mu$-stable vector bundles $E[1]$. Note
that a torsion free sheaf is $\mu$-stable if and only if it has
no proper subsheaf. (The degree of all sheaves is trivial!)

Stable objects and minimal objects coincide for $z\in R_0$, but this
does not hold for $z\in R_\pm$.

iii) Later in Lemma \ref{lem:autogen} it will be important that
for $z\in R_0$ the only stable semi-rigid objects are the point
sheaves $k(x)$. Here one uses that on a generic K3 surface any
stable vector bundle $E$ of rank $>1$ has $(E,E)^1\geq4$.

The same assertion does not hold for $z\in R_\pm$.
Indeed, for $z\in R_+$ the
objects given by a non-split exact sequence $0\to \ko_X[1]\to F\to
k(x)\to 0$ are semi-rigid and stable in $\ka$
and the ideal sheaves $I_x\in \coh(X)$
of points $x\in X$ are semi-rigid and stable in $\coh(X)$.

%
% all objects in
%$\ka(x,y=0)$ have phase zero and are, therefore, semistable. The
%stable ones are the simple objects of the category, i.e.\ the
%point sheaves $k(x)$ and the shifts $E[1]$ of $\mu$-stable vector
%bundles $E$.
%%
%
%In i) again all the point sheaves $k(x)$, $x\in X$ and all $E[1]$
%with $E$ an arbitrary stable vector bundle are stable objects in
%$\ka(x,y)$, although this time $\phi(E[1])<1$. They form once
%more the simple objects of the category. Note however, that many
%more stable objects can be produced by interpreting a non-trivial
%class in $\Ext^2(k(x),E)$ as an extension $0\to E[1]\to F\to
%k(x)\to 0$ (here $E$ and $k(x)$ are as before) and $F$ is as well
%stable.
%
%In the last case iii), where $\ka(x,y)=\coh(X)$, the stable
%objects are the point sheaves $k(x)$ and the unshifted(!) stable
%vector bundles.
%
%In all three cases, the only minimal%
%
%From here one finds that in the cases ii) and iii), i.e.\
%$y\leq0$, the only stable semi-rigid objects in $\ka(x,y)$ are
%the point sheaves $k(x)$. Indeed,
%
%
\end{remark}

\begin{prop} For $z\in R$ the stability function $Z:\ka(z)\to \CC$ has the
Harder--Narasimhan property. Moreover, the induced stability condition
is locally finite and in particular all Jordan--H\"older
filtrations are finite.
\end{prop}

\begin{proof}
For $z\in R_-$, this follows from the existence of the usual
Harder--Narasimhan  filtration in $\ka(z)=\coh(X)$. In fact, as
$\Pic(X)=0$ the Harder--Narasimhan filtration consists of the
torsion part and the rest. Furthermore, due to the special form
of the stability function one can easily verify that the stability
condition is locally finite.

For $z\in R_+\cup R_0$ one easily proves that (independently of any
stability condition) the abelian category $\ka=\ka(z)$ is of finite
length, i.e.\ any object admits a finite filtration
whose quotients are minimal objects of the category.
This follows from the description of the minimal objects
given above.
Consequently, any filtration of an object in $\ka$
becomes stationary. This applies in particular to
Harder--Narasimhan and Jordan--H\"older filtrations.
\end{proof}

So, we can conclude that for any $z\in R$ we have defined a
locally finite stability condition $$\sigma_{z}\in\Stab(\Db(X))$$
given by the $t$-structure associated to the torsion theory
$(\kf(z),\kt(z))$ and the stability function $Z$ on its heart.

This allows us to consider $R$ in a natural way as a subset of
the space of stability conditions:$$R\subset \Stab(X).$$ For any
$z\in R$ the point sheaves $k(x)$ are $\sigma_{z}$-stable of phase
one with $Z(k(x))=-1$ and if, moreover, $z\in R_0$, then they are
the only semi-rigid stable objects (up to shift).

\begin{remark}
Note that due to \cite{B}, $\coh(X)$ never occurs as the heart of
a stability condition on a projective K3 surface. Indeed,  if
$\coh(X)=\kp((0,1])$, then the point sheaves $k(x)$, which are
minimal objects in $\coh(X)$, would all be stable of the same
maximal phase $\phi\in(0,1]$. Then one easily shows $\phi(k(x))=1$
and \cite[Prop.\ 10.3]{B} yields $\kp((0,1])=\ka(\varphi)$ for some
$\varphi=\exp(B+i\omega)$. Clearly, if $X$ is projective, then
$\ka(\varphi)\ne\coh(X)$.
\end{remark}

\bigskip
%%%%%%%%%%%%%%%%%%%%%%%%%%%%%%%%%%%%%%%%%%%%%%%%%%%%%%%%%%%%%%%%%%%%%

\subsection{The space of all stability conditions on
a generic K3 surface}\label{sect:Class}~

\smallskip

In the next step, again following Bridgeland, one tries to
classify all stability conditions $\sigma=(Z,\kp)$ for which all
point sheaves $k(x)$ are stable of phase one. Up to scaling we
may assume that $Z(k(x))=-1$. Hence the stability function $Z$ is
of the form considered above: $Z(E)=-u\cdot r-s-i(v\cdot r)$.

By Proposition \ref{prop:stabonly} the trivial line bundle
$\ko_X$, which is up to shift the only spherical object in
$\Db(X)$, is stable with respect to $\sigma$. Since
$(\ko_X,k(x))^0\ne0\ne(k(x),\ko_X)^2$, one has
$\phi(\ko_X)<1<\phi(\ko_X)+2$. Hence, either $\ko_X\in\kp((0,1])$
or $\ko_X[1]\in\kp((0,1])$ depending on whether $v<0$ or $v\geq0$.
If $v=0$, then necessarily $u<-1$. Note that so far we have only
used that there is one point sheaf $k(x)$ which is
$\sigma$-stable of phase one.

Now one copies the proof of Lemma 10.1 and Proposition 10.3 in \cite{B}.
In fact the arguments simplify, as $\Pic(X)=0$. In particular, the
proof of \cite[Prop.\ 10.3]{B} does not only show equality of the
hearts, but right away equality of the stability conditions.

\begin{prop}\label{prop:classif}
Suppose $\sigma=(Z,\kp)$ is a stability condition on $\Db(X)$ for
a K3 surface $X$ with trivial Picard group. If all point sheaves
$k(x)$ are stable of phase one with $Z(k(x))=-1$, then
$\sigma=\sigma_{z}$ as above for some $z\in R$.\qqed
\end{prop}

For later use we point out that the heart $\kp((0,1])$ of
$\sigma$ is $\ka$ if ${\rm Im}(Z(\ko_X))\leq0$ and $\coh(X)$ if
${\rm Im}(Z(\ko))>0$.

Combining Proposition \ref{prop:classif}
 with Corollary \ref{cor:Tnstableonly} we get a
complete classification of all stability conditions on $\Db(X)$
for a generic K3 surface $X$. As before, $T$ denotes the spherical
shift $T_{\ko_X}$.

\begin{cor}\label{CorStab}
Suppose $\sigma$ is a stability condition on $\Db(X)$ for a K3
surface $X$ with trivial Picard group. Then up to the action of
$\widetilde{\rm Gl}_2^+(\RR)$ the stability condition $\sigma$ is
of the form $T^n(\sigma_{z})$ for some $z\in R$ and
$n\in\ZZ$.\qqed
\end{cor}

\begin{thm}\label{thm:classgen}
If $X$ is a K3 surface with trivial Picard group, then the space $\Stab(X)$
of stability conditions on $\Db(X)$ is connected and simply-connected.
\end{thm}

\begin{proof} Let us consider $$W(X):=\widetilde{\rm Gl}_2^+(\RR)(R)\subset \Stab(X),$$
which can also be written as the union $W(X)=W_+\cup W_-\cup W_0$ according
to (\ref{decompR}). Then Corollary \ref{CorStab} says
\begin{equation}\label{opencover}\Stab(X)=\bigcup_n T^n W(X).
\end{equation}

Recall that the group $\widetilde{\rm Gl}_2^+(\RR)$
can be thought of as the set of pairs $\tilde g=(g,f)$ of a linear map
$g\in {\rm Gl}_2^+(\RR)$ and an increasing map $f:\RR\to \RR$ with
$f(\phi+1)=f(\phi)+1$ inducing the same map on
${\rm S}^1=\RR/2\ZZ=\RR^2/\RR_{>0}$. Using the natural action
of ${\rm Gl}_2^+(\RR)$ on $\CC$, one defines the action of an element $\tilde g$
on $\Stab(X)$ as follows: If
$\sigma=(Z,\kp)$, then $\tilde g(\sigma)=(Z',\kp')$ satisfies $Z'(E)=g^{-1}Z(E)$
and $\kp'(\phi)=\kp(f(\phi))$.
The action is clearly continuous. See \cite[Lemma 8.2]{B2}.

i) Suppose that for some $\sigma_z\in R$ and $\tilde g\in
\widetilde{\rm Gl}_2^+(\RR)$ one has $\tilde g(\sigma_z)\in R$.
Then $g^{-1}=\left(\begin{array}{cc}1&a\\0&b\end{array}\right)$
with $b>0$. Thus, if $\sigma_{z'}=\tilde g(\sigma_z)$, then
$z'=u+av+ibv$. In particular, the sign of the imaginary part does
not change. Hence, $W(X)=W_+\cup W_-\cup W_0$ is disjoint.

The calculation also shows that $W_+$ and $W_-$ are  contained in
two different orbits of the $\widetilde{\rm Gl}_2^+(\RR)$-action,
whereas two points in $R_0$ are contained in two different orbits.
More precisely, if $z\in R_\pm$, then the continuous surjection
$\{g\in{\rm Gl}_2^+(\RR)~|~g^{-1}=
\left(\begin{array}{cc}1&a\\0&b\end{array}\right)\}\twoheadrightarrow
R_\pm$, $g\mapsto g^{-1}z$ also describes the inclusion
$R_\pm\subset \Stab(X)$ by lifting $g$ to  $\tilde g=(g,f)$ such
that $f((0,1])=(0,1]$ and mapping $\sigma_z$ to $\tilde g(
\sigma_z)=\sigma_{g^{-1}z}$.

\smallskip

ii) We claim that $W(X)\subset\Stab(X)$ is an open connected subset.

Let us first prove that the inclusion $R\subset\Stab(X)$ is
continuous. The arguments in i) show that the restrictions
$R_\pm\to\Stab(X)$ are continuous. Thus, it suffices to prove that the inclusion
is continuous in points $z\in R_0$.

For this consider a $\sigma_z$-stable object $0\ne F$ of phase $\phi$. Hence
$\phi=\ell$ for some $\ell\in\ZZ$ and $F\cong k(x)[\ell-1]$ with
$x\in X$ or $F=E[\ell]$ with $E$  stable and locally free (see
Remark \ref{rem:onlysemirig}, ii)). Obviously both types of
objects stay $\sigma_{z'}$-semistable for $z'\in R_\pm$ close to $z$
with phase $\phi'$ close to $\phi=\ell$. The important point to notice here is that
how close $\phi'$ is to $\phi$ only depends on $z'$ and not on $F$. This is
clear for $F\cong k(x)[\ell-1]$, whose phase stays $\ell$. For $F\cong E[\ell]$
this is a priori not so clear. In order to see this, one has
to check that  $vr/(ru+s)$ is uniformly small (independent of $r=\rk(E)$
and $s=\rk(E)- {\rm c}_2(E)$)  for $|u|<\delta <\!\!<0$. This follows e.g.\
from ${\rm c}_2(E)\geq2\rk(E)$ (see the proof of Lemma \ref{lem:stabfctgen}),
as $vr/(ru+s)=v/(u+1-{\rm c}_2(E)/\rk(E))$.

Thus, $R$ and hence $W(X)=\widetilde{\rm Gl}_2^+(\RR)(R)$ are
connected subsets of $\Stab(X)$. Moreover,
$W(X)\to\kn(X)^*\otimes \CC\cong\CC^2$ is a local homeomorphism.
Since the same is true for  the period map $\Stab(X)\to
\kn(X)^*\otimes\CC$ (see \cite[Thm.\ 1.2]{B2}), this suffices to
conclude the openness of $W(X)$ inside $\Stab(X)$ from the
openness of its image in $\CC^2$.

\smallskip

iii) Due to Remark \ref{rem:t2not}, ii) one knows that $T^nW(X)$
and $T^kW(X)$ are disjoint for $|n-k|\geq2$. We now claim that
\begin{equation}\label{interequ}
T^nW(X)\cap T^{n+1}W(X)=T^nW_-.
\end{equation}

Of course, it is enough to prove $W(X)\cap TW(X)=W_-$.  More precisely, we shall
show $$TW_+=W_-~~\text{and}~~T(W_0\cup W_-)\cap W(X)=\emptyset.$$

Let us consider $\sigma:=\sigma_i=(Z,\kp)\in R_+$ and $\sigma':=\sigma_{-i}=(Z',\kp')\in R_-$.
We shall show that $T\sigma=\tilde i\cdot\sigma'$, where $\tilde i\in\widetilde{\rm Gl}_2^+(\RR)$
acts by $\tilde i(Z')(E)=-i Z'(E)$ and $(\tilde i\kp')(\phi)=\kp'(\phi+1/2)$.

The verification of $T(Z)(E)=-iZ'(E)$ can be done easily by using
that the spherical twist acts on cohomology by $(r,0,s)\mapsto
-(s,0,r)$. For the second equality it will be enough to show that
$T\kp((0,1])=(\tilde i\kp')((0,1])=\kp'((1/2,3/2])$. Clearly, two
stability conditions with identical  hearts and stability
functions coincide.

By construction $\ka=\kp((0,1])$ and $\coh(X)=\kp'((0,1])$.
Moreover, for any closed point $x\in X$ one has
$T(k(x))=I_x[1]\in\kp'(1/2)[1]=\kp'(3/2)$. Thus,
the heart $T^{-1}\kp'((1/2,3/2])$ of the stability condition
$T^{-1}(\tilde i\sigma')$ contains
all point sheaves $k(x)$ which, moreover, have central charge $-1$.
By Proposition \ref{prop:classif} this is enough to conclude
$T^{-1}\kp'((1/2,3/2])=\ka$. (A priori
one could have $T^{-1}\kp'((1/2,3/2])=\coh(X)$, but this is impossible as $T\ko=\ko[-1]\in
\kp'(-1/4)\not\subset\kp'((1/2,3/2])$.)

Next, $TW_-=T^2W_+$, but due to  Remark \ref{rem:t2not}, ii) one
has $T^2W_+\cap W(X)=\emptyset$. Eventually, $TW_0\cap
W(X)=\emptyset$, because the stability function of $T\sigma_z$
with $z=u<-1$ is given by $(r,0,s)\mapsto u\cdot s+r$, which is
not of the type realized by any stability condition in $W(X)$.

\smallskip

iv) As an immediate consequence of (\ref{interequ}) in iii) and
the connectedness of $W(X)$ proved in ii),
 one concludes that $\Stab(X)=\bigcup T^nW(X)$ is connected.

v) In the final step one applies the van Kampen theorem to the open cover
(\ref{opencover}). The intersections $T^nW(X)\cap T^{k}W(X)\subset\Stab(X)$ are
either empty for $|n-k|\geq2$ or homeomorphic to the connected $W_-$.
Thus, it suffices to verify that the open sets $T^nW(X)\cong W(X)$
are simply-connected.

For this purpose  let us consider its image in $\CC^2$. Clearly,
${\rm Gl}_2^+(\RR)(R_+\cup R_-)\subset\CC^2$ is the set of all
$\RR$-linearly independent pairs $(z_0,z_1)\in\CC^2$ and ${\rm
Gl}_2^+(\RR)(R_0)$ consists of pairs $(z_0,uz_0)$ with $z\in\CC^*$
and $u\in(-\infty,-1)$. Thus, ${\rm Gl}_2^+(\RR)(R)$ is the
complement of
$\{(z,uz)~|~z\in\CC,u\in[-1,\infty)\}\bigcup\{(0,z)~|~z\in\CC\}$,
whose fundamental group is generated by the loop around the real
codimension-one component $\{(0,z)\}$. This loop can be written as
the rotation $g_t$ by $\pi t$, $t\in[0,2]$ which lifts to
$\tilde g_t=(g_t,f_t)\in\widetilde {\rm Gl}_2^+(\RR)(R)$ with
$f_t(\phi)=\phi+t$. In particular, $\widetilde g_2$ is the double
shift acting as a deck transformation on $W(X)$. Thus, $W(X)$ is
simply-connected.
\end{proof}

Some of the arguments in the above proof are inspired by a very detailed description of
all stability conditions on a generic K3 surface due to S.\ Meinhardt \cite{M}. In analogy to \cite{B},
he introduces a period domain and presents $\Stab(X)$ via a period
map as the universal cover of it. The main difference compared to \cite{B}
is that the naive definition of the period domains
$P(X)$ or $P^+(X)$ makes no sense, as there
are no positive planes in $\kn(X)\otimes\RR$.

\bigskip
%%%%%%%%%%%%%%%%%%%%%%%%%%%%%%%%%%%%%%%%%%%%%%%%%%%%%%%%%%%%%%%%%%%%%

\subsection{Autoequivalences}\label{subsect:autononproj} ~

\smallskip

A complete description of the group of autoequivalences of
$\Db(X)$ in the generic case can now either be obtained by
following the general approach in \cite{B} or by a more direct
argument given below, which does not rely on the description of
$\Stab(X)$ in Theorem \ref{thm:classgen}.

\begin{lem}\label{lem:autogen}
If $\Phi_\ke:\Db(X)\congpf\Db(Y)$ is a Fourier--Mukai equivalence
between two K3 surfaces $X$ and $Y$ with $\Pic(X)=0$, then  up to
shift
$$\Phi_\ke\cong T^n\circ f_*$$
for some $n\in\ZZ$. Here, $T$ is the spherical twist with respect
to $\ko_Y$ and $f:X\congpf Y$ is an isomorphism.
\end{lem}

\begin{proof}
Suppose $\sigma=\sigma_{(u,v=0)}$ is one of the distinguished
stability conditions constructed above with $(u,v)\in R$ and
$v=0$. Let $\tilde\sigma$ be its image under $\Phi_\ke$.

Then by Proposition \ref{prop:onlyone} and Corollary
\ref{cor:Tnstableonly} there exists an integer $n$ such that all
point sheaves $k(x)$ are stable of the same phase with respect to
$T^n(\tilde\sigma)$. Thus the composition
$\Psi:=T^n\circ\Phi_\ke$, which is again of Fourier--Mukai type,
i.e.\ $\Psi=\Psi_\kf$ for some $\kf\in\Db(X\times Y)$, sends the
stability condition $\sigma$ to a stability condition $\sigma'$
for which all point sheaves are stable of the same phase. Shifting
the kernel $\kf$ allows one to assume that
$\phi_{\sigma'}(k(y))\in(0,1]$ for all $y\in Y$. Thus, the heart
$\kp((0,1])$ of $\sigma'$, which under $\Psi_\kf$ is identified
with $\ka(u)$, contains as stable objects the images $\Psi(k(x))$
of all points $x\in X$ and all point sheaves $k(y)$.

But as was remarked earlier (see Remark \ref{rem:onlysemirig}, iii),
the only semi-rigid stable objects in $\ka(u)$ are the point
sheaves. Thus, $\Psi^{-1}(k(y))$ must be of the form $k(x)$. In
other words, for all $y\in Y$ there exists a point $x\in X$ such
that $\Psi(k(x))\cong k(y)$.

This suffices to conclude that the Fourier--Mukai equivalence
$\Psi_\kf$ is a composition of $f_*$ for some isomorphism
$f:X\congpf Y$ and a line bundle twist (cf.\ \cite[Cor.\
5.23]{HFM}), but there are no non-trivial line bundles on $Y$.
\end{proof}

This immediately leads to the following complete description of
all Fourier--Mukai equivalences in the generic case.

\begin{prop}\label{prop:autgennonproj}
If $\Aut(\Db(X))$ denotes the group of autoequivalences of
Fourier--Mukai type of a K3 surface $X$ with $\Pic(X)=0$, then
$$\Aut(\Db(X))\cong\ZZ\oplus\ZZ\oplus\Aut(X).$$ The first two
factors are generated respectively by the shift functor and  the
spherical twist $T_{\ko_X}$.\qqed
\end{prop}

Maybe it is worth pointing out that in the non-algebraic case not
every derived equivalence is of Fourier--Mukai type (see
\cite{V4}). Also, the automorphism group of a K3 surface with
trivial Picard group can be explicitly described. Is either
trivial or isomorphic to $\ZZ$. See \cite[Cor.\ 1.6]{Og}.
\bigskip

%%%%%%%%%%%%%%%%%%%%%%%%%%%%%%%%%%%%%%%%%%%%%%%%%%%%%%%%%%%%%%%%%%%%%%%%%

\end{document}